\documentclass[12pt]{article}

\usepackage{amsmath}
\usepackage{amsfonts}
\usepackage{amsthm}
\usepackage{hyperref}

\begin{document}

\newcommand{\A}{\mbox{${{{\cal A}}}$}}
\newcommand{\Ato}{\mbox{ ${{{\overset{\A}{\longrightarrow}}}}$ }}

%1111111111111111111111111111111111111111111111111111111111111111111111111111

\author{Attila Losonczi}
\title{A common generalization of infinite sum, unordered sum and integral}

\date{\today}%{8 September 2021}

\newtheorem{thm}{\qquad Theorem}[section]
\newtheorem{prp}[thm]{\qquad Proposition}
\newtheorem{lem}[thm]{\qquad Lemma}
\newtheorem{cor}[thm]{\qquad Corollary}
\newtheorem{rem}[thm]{\qquad Remark}
\newtheorem{ex}[thm]{\qquad Example}
\newtheorem{df}[thm]{\qquad Definition}
\newtheorem{prb}{\qquad Problem}

\newtheorem*{thm2}{\qquad Theorem}
\newtheorem*{cor2}{\qquad Corollary}
\newtheorem*{prp2}{\qquad Proposition}

\maketitle

\begin{abstract}

\noindent

We present a common ground for infinite sums, unordered sums, Riemann/Lebesgue integrals, arc length and some generalized means. It is based on extending functions on finite sets using Hausdorff metric in a natural way.

\noindent
\footnotetext{\noindent
AMS (2020) Subject Classifications: 30L99, 28A10, 40A05, 26E60 \\

}

\end{abstract}

%-------------------------------------------------------------------------------------------------------------------------------------------Introduction----------------------------------
\section{Introduction}
In Analysis there are many constructions which somehow extend finite summation, finite average or finite weighted average in a natural way. Let us enumerate a few examples.

\begin{itemize}

\item The sum of an infinite series consisting of positive elements is an extension of summation on finite (multi-) sets. Or more generally the sum of an absolute convergent series is the same type of extension. In these cases the sum does not depend on the order hence we can talk about the sum of the elements of the series without ordering them. 

\item If we have a set $H\subset\mathbb{R}$ then we can define the unordered sum of the elements in $H$ or more generally the unordered sum of a function can be taken (see section \ref{sbnn}). It is again a kind of extension of summation on finite sets.

\item Similarly to the previous notion we can define the unordered average of the elements in a set $H$ or more generally the unordered average of a function can be taken (defined and examined later in this paper in section \ref{unorderedavg}). It is a kind of extension of average on finite sets.

\item The general concept of the sum of infinite series extends finite summation in a slightly different way: here the order does matter. I.e. it matters in which order we estimate the infinite sum by finite sums, in which order we take finite subsets of the elements.       

\item The Riemann integral of a function is a kind of extension of weighted averages of function values on finite sets (on intervals with length equals to 1).

\item The integral of a function on a measure space is a kind of extension of special weighted averages of values on finite sets (on spaces with measure equals to 1).

\item The arc length of a continuous curve can be defined in a way too such that it is a kind of extension of some values on finite sets (see section \ref{sbnn}).

\item Let $H\subset\mathbb{R}$ be bounded. If the isolated points determine the set in the sense that $cl(H-H')=H$ then the generalized mean ${\cal{M}}^{iso}(H)$ is calculated by the limit of averages on certain finite subsets of the isolated points (see definition in \ref{defiso}).

\item  Let $H\subset\mathbb{R}$ be bounded. Then the generalized mean ${\cal{M}}^{eds}(H)$ is defined as the limit of averages of special evenly distributed finite samples from $H$ (see definition in \ref{defeds}).

\end{itemize}

The natural question arises that can one find a common generalization of these notions? The answer is yes and based on extending functions on finite sets which approximate the underlying set in Hausdorff metric. 

\smallskip

After some introduction in the first part of the paper we define and examine some basic notions e.g unordered average.

In the second part of the paper we define the general extension method then we show that is applicable in all mentioned cases.

In the third part of the paper we establish some general properties of that extension method. For example we examine uniqueness, relations to subspaces, closer of sets, extensions to smaller sets and among others, continuity. We also present many examples.

At the end of the paper we provide some generalization of the already developed notions in order to be able to handle integral on measure space in general. 

%------------------------------------------------------------------------------------------------------------------------------Basic notions and notations---------------------------
\subsection{Basic notions and notations}\label{sbnn}

In this subsection we enumerate many (well) known notions for the paper to be more or less self-contained.

\medskip

A \textbf{multiset} (or mset) is a pair $\langle H,m\rangle$ where $H$ is a set and $m:H\to\mathbb{N}$. It generalizes the common notion of set in a way that it allows multiple occurrences of elements of a set (see \cite{syr},\cite{blizard}). We say that $\langle K,l\rangle$ is a subset of $\langle H,m\rangle$ if $K\subset H$ and for each $k\in K\ \ l(k)\leq m(k)$ holds. For multisets we use the notation $\langle K,l\rangle\subset \langle H,m\rangle$.

\medskip

The following definition of \textbf{unordered sum} is from \cite[Definition 3.2]{bruckner} but is included here for easier readability.
Let $a:I\to\mathbb{R}$. Then we write $\sum\limits_{i\in I}a(i)=c$ if $\forall\epsilon>0$ there exists a finite set $I_0\subset I$ such that for every finite set $J$, $I_0\subset J\subset I$, $\left|\sum\limits_{i\in J}a(i)-c\right|<\epsilon$ holds.

It is known (see \cite{bruckner}) that if the unordered sum converges then only countably many of its terms can be non zero, moreover the only accumulation point of its terms is 0. Moreover its terms can be arranged in a series that is absolute convergent with exactly the same sum. It is also known that the sum of an absolute convergent series is not sensitive for reordering therefore it is readily seen that the following definition gives equivalent notion to both the unordered sum and the sum of an absolute convergent series.

\begin{df}Let $\langle H,m\rangle$ be a multiset. Then we write $\sum\limits_{h\in H}h=c$ if $\forall\epsilon>0$ there exists a finite mset $\langle H_0,m_0\rangle\subset \langle H,m\rangle$ such that for every finite mset $\langle J,m_j\rangle$, $\langle H_0,m_0\rangle\subset \langle J,m_j\rangle\subset \langle H,m\rangle$ implies that $\left|\sum\limits_{h\in J}m_j(h)\cdot h-c\right|<\epsilon$.
\end{df}

Throughout this paper function $\A(\dots)$ will denote the arithmetic mean of any (finite) number of variables.  If $K\subset\mathbb{R}$ finite then $\A(K)$ will denote the arithmetic mean of the elements of $K$ (the order obviously does not matter).

\smallskip

$cl(H), H'$ will denote the closure and accumulation points of a subset $H$ of $\mathbb{R}$ respectively. 

\smallskip

We use the usual notation: $\overline{\mathbb{R}}=\mathbb{R}\cup\{-\infty,\infty\}$.

\smallskip

In the sequel let $\langle I,d\rangle$ denote a pseudo-metric space. We use the usual notations:
\[d(x,H)=\inf\{d(x,y):y\in H\},d(K,H)=\inf\{d(x,y):x\in K,y\in H\}\ (K,H\subset I),\]
\[S(x,\epsilon)=\{y\in I:d(x,y)<\epsilon\},\ S(A,\epsilon)=\{y\in I:d(A,y)<\epsilon\}=\bigcup\limits_{x\in A}S(x,\epsilon).\]
\[B(x,\epsilon)=\{y\in I:d(x,y)\leq\epsilon\},\ B(A,\epsilon)=\{y\in I:d(A,y)\leq\epsilon\}.\]
If we want to empathize that $S$ or $B$ is taken in the space $I$ then we use the notations $S_I$, $B_I$ respectively.

We recall the definition of \textbf{Hausdorff pseudo-metric} on the subsets of the pseudo-metric space $\langle I,d\rangle$: Let $A,C\subset I$. Then 
\[d_H(A,C)=\inf\big\{\epsilon>0: A\subset S(C,\epsilon),\ C\subset S(A,\epsilon)\big\}\]
is a pseudo-metric (see \cite{ert}). It is easy to see that $d_H(A,C)=\inf\{\epsilon>0: A\subset B(C,\epsilon),\ C\subset B(A,\epsilon)\}$.

\medskip

The following definitions are from \cite{lamis} but is included here for easier readability.

\begin{df}\label{defiso}Let $H\subset\mathbb{R}$ be bounded and $cl(H-H')=H$. Then let \[{\cal{M}}^{iso}(H)=\lim_{\delta\to 0+0}\A\big(H-S(H',\delta)\big)\] if the limit exists. 
\end{df}

Note that $H-S(H',\delta)$ is finite for all $\delta>0$.

\begin{df}\label{defeds}Let $H\subset\mathbb{R}$ be bounded, $a=\inf H,b=\sup H$. If $n\in\mathbb{N},0\leq i\leq n$ then set $H_{n,i}=H\cap\big[a+\frac{i}{n}(b-a),a+\frac{i+1}{n}(b-a)\big)$. Let $\ I_n=\{0\leq i\leq n:H_{n,i}\ne\emptyset\}$.

We say that the mean of $H$ is $k={\cal{M}}^{eds}(H)$ if $\forall\epsilon>0\ \exists N\in\mathbb{N}$ such that $n>N, \xi_i\in H_{n,i}\ (i\in I_n)$ implies that $\big|\A\big(\{\xi_i:i\in I_n\}\big)-k\big|<\epsilon$.
\end{df}

${\cal{M}}^{iso}(H)$ and ${\cal{M}}^{eds}(H)$ are generalized means of the set $H$ in the sense of \cite{lamis}.

\medskip

If $I$ is a set then set \[[I]^{<\infty}=\{H\subset I: H\text{ is finite}\},\ [I]^{n}=\{H\subset I: |H|=n\}\ (n\in\mathbb{N}).\]
If $s:[I]^{<\infty}\to\mathbb{R}$ and $J\subset I$ then let us use the notation $s||_J=s|_{[J]^{<\infty}}$ for the restriction of $s$ for the finite subsets of $J$.

\smallskip

Let $\gamma:[0,1]\to\mathbb{R}^n$ be continuous. The usual definition for the \textbf{arc length} of $\gamma$ is the supremum of the length of all finite polygons inscribed in $\gamma$. An alternative definition can be the following. Let $I=[0,1]$. The arc length of $\gamma$ is $A\in\mathbb{R}$ if $\forall\epsilon>0\ \exists\delta>0$ such that $K\in[I]^{<\infty}$ and $d_H(K,I)<\delta$ implies that the distance between $A$ and the length of the inscribed polygon determined by $K$ is less than $\epsilon$. (Similar way works when $A=\pm\infty$.) We will show that the two definitions are equivalent (see \ref{carclendefeq}).

\smallskip

Let $\langle I_i,d_i\rangle$ be pseudo-metric spaces for $i=1,2$. Then it is known that $\langle I_1\times I_2,d\rangle$ is the \textbf{product space} for $d\big((x_1,y_1),(x_2,y_2)\big)=d_1(x_1,x_2)+d_2(y_1,y_2)\ \ (x_1,x_2\in I_1,\ y_1,y_2\in I_2)$. In the sequel we will always consider that metric on the product.

\smallskip

If $x\in\mathbb{R}$ then we use the usual notations $|x|^+=\max\{0,x\},\newline |x|^-=-\min\{0,x\}$.

\smallskip

If $a\in\mathbb{R}$ then let $a+(+\infty)=+\infty$, if $a\in\mathbb{R}^+$ let $a\cdot (+\infty)=+\infty$, and let $0\cdot (+\infty)=0$.  

%------------------------------------------------------------------------------------------------------------------------------On $\delta$-dense sets---------------------------
\section{On $\boldsymbol\delta$-dense sets}

In the sequel $\langle I,d\rangle$ (or simply $I$) denotes a pseudo-metric space.

\smallskip

A subset $K$ of $I$ is $\boldsymbol\delta$\textbf{-dense} (for $\delta>0$) if $\forall x\in I\ \exists k\in K$ such that $d(k,x)<\delta$ or equivalently $S(K,\delta)=I$  (see \cite{ert}).

In a bounded space every subset is $\delta$-dense for some $\delta>0$.

\smallskip

$I$ is \textbf{totally bounded} if $\forall \delta>0$ there exists a finite subset $K$ of $I$ such that $K$ is $\delta$-dense (see \cite{ert}).

\begin{prp}If $K$ is $\delta$-dense then $d_H(K,I)\leq\delta$. 

If $K$ is finite and $d_H(K,I)\leq\delta$ then $\forall x\in I\ \exists k\in K$ such that $d(k,x)\leq\delta$.
\end{prp}
\begin{proof}The first statement is straightforward. 

To prove the second one, suppose that there is $x\in I$ such that $\forall k\in K\ d(k,x)>\delta$. Let $\delta_1=\min\{d(k,x):k\in K\}$. Then $S(K,\delta_1)\ne I$ which gives that $d_H(K,I)\geq\delta_1$ -- a contradiction.
\end{proof}

\begin{ex}Let $I=\big\{\frac{1}{n}:n\in\mathbb{N}\big\}\cup\{-1\}$ and $K=\big\{\frac{1}{n}:n\in\mathbb{N}\big\}$. Then $d_H(K,I)=1$ and $K$ is not $1$-dense.
\end{ex}

\begin{df}A subset $K$ of $I$ is called $\boldsymbol\delta$\textbf{-strong-dense} or simply $\boldsymbol\delta$\textbf{-sdense} if there is $\delta'<\delta$ such that $K$ is $\delta'$-dense.
\end{df}

\begin{prp}Let $K\subset I$. Then $K$ is $\delta$-sdense iff $d_H(K,I)<\delta$.
\end{prp}
\begin{proof}If $K$ is $\delta$-sdense then there is $\delta'<\delta$ such that $K$ is $\delta'$-dense i.e. $S(K,\delta')=I$ which gives that $d_H(K,I)\leq\delta'$.

If $d_H(K,I)<\delta$ then there is $\delta'<\delta$ such that $d_H(K,I)=\delta'$. If $\delta_1=\frac{\delta+\delta'}{2}$ then $S(K,\delta_1)=I$ hence $K$ is $\delta_1$-dense. 
\end{proof}

\begin{df}In a bounded space $I$, $K\subset I$ is \textbf{stretched} if for all $x,y\in K,\ x\ne y,\ d(x,y)\geq d_H(K,I)$ holds. 

$K\subset I$ is called \textbf{strongly stretched} if $x,y\in K,\ x\ne y,\ d(x,y)>d_H(K,I)$ holds.

Let $[I]_{str}$ and $[I]_{sstr}$ denote all stretched and strongly stretched subsets of $I$ respectively and let $[I]^{<\infty}_{str}=[I]_{str}\cap[I]^{<\infty},\ [I]^{<\infty}_{sstr}=[I]_{sstr}\cap[I]^{<\infty}$.\qed
\end{df}

\begin{prp}\label{pexistsstr}If $I$ is totally bounded then there exists a finite, strongly stretched, $\delta$-dense set for each $\delta>0$.
\end{prp}
\begin{proof}Let $K$ be a finite $\frac{\delta}{4}$-dense set in $I$. Let $p_1\in I$ be arbitrary. If $p_1, \dots, p_n$ have already been chosen then let $p_{n+1}\in I-\bigcup\limits_{i=1}^nB(p_i,\frac{\delta}{2})$. This process has to finish in at most $|K|$ many steps because if $p,q\in S(h,\frac{\delta}{4})$ then $d(p,q)<\frac{\delta}{2}$ therefore $S(h,\frac{\delta}{4})$ can contain at most 1 point of $\{p_1, \dots, p_n\}$. Hence $\bigcup\limits_{i=1}^nB(p_i,\frac{\delta}{2})=I$ for some $n$. Clearly $d_H(\{p_1, \dots, p_n\},I)\leq\frac{\delta}{2}$, $S(\{p_1, \dots, p_n\},\delta)=I$ and $d(p_i,p_j)>\frac{\delta}{2}\ (i\ne j)$ therefore $\{p_1, \dots, p_n\}$ is a strongly stretched $\delta$-dense set.
\end{proof}

\begin{prp}If $I$ is totally bounded $K\in[I]^{<\infty}$ then there is $L\in[I]^{<\infty}_{str}$ such that $K\subset L$.
\end{prp}
\begin{proof}Let $\delta=\min\{d(x,y):x,y\in K,x\ne y\}$. If $\delta\geq d_H(K,I)$ then $K$ is already stretched. Suppose that $\delta<d_H(K,I)$. Let $K=\{p_1, \dots, p_m\}$ as starting points and then follow the method in \ref{pexistsstr} and then we will end up with a stretched $\delta$-dense set that contains $K$.
\end{proof}

\begin{prp}If $I$ is a bounded pseudo-metric space then 
\[\sup\{|K|:K\in[I]_{str},\ d_H(K,I)=2\delta\}\leq\min\{|K|:K\subset I,\ d_H(K,I)=\delta\}.\]
\end{prp}
\begin{proof}Take a $K\in[I]_{str},\ d_H(K,I)=2\delta$ and $L\subset I,\ d_H(L,I)=\delta$. If $k\in K$ then there is a $l_k\in L$ such that $k\in S(l_k,\delta)$. Clearly $k\ne m\ (k,m\in K)$ implies that $l_k\ne l_m$ which gives that $|K|\leq|L|$.
\end{proof}

\begin{cor}If $I$ is a totally bounded pseudo-metric space then 
\[\max\{|K|:K\in[I]^{<\infty}_{str},\ d_H(K,I)=2\delta\}\leq\min\{|K|:K\in[I]^{<\infty},\ d_H(K,I)=\delta\}.\]
\end{cor}

\begin{cor}For given $\delta$ there is a $K\in[I]^{<\infty}_{str},\ d_H(K,I)=\delta$ with maximal cardinality.\qed
\end{cor}

\begin{cor}If $I$ is a totally bounded pseudo-metric space and $K$ is stretched then $K$ is finite.\qed
\end{cor}

%\begin{prp}If there is a finite $\delta$-dense set then there exists a finite, stretched $\delta$-dense set as well.
%\end{prp}
%\begin{proof}Let $H$ be a finite $\delta$-dense set in $I$. ????
%\end{proof}

\begin{ex}It is not true that a finite $\delta$-dense set has a stretched $\delta$-dense subset.  

Let $I=[0,1], H=\{0.25,0.75\}$. Then $H$ is $0.5$-dense but there is no stretched $0.5$-dense subset of $H$.
\end{ex}

\begin{ex}It is not true that two finite, stretched and $\delta$-dense sets have the same cardinality.

Let $I=[0,1], H_1=\{0,0.5,1\}, H_2=\{0.25,0.75\}$. Then both sets are stretched and $\delta$-dense for all $\delta>0.25$.
\end{ex}

\begin{ex}Let $\delta=d(K,I)$ and 
\[|K|=\min\{|L|: L\in[I]^{<\infty},\ \delta\text{-dense}\}.\] 
It does not follow that $K$ has to be stretched.

Let $I=(0,1.4), \delta=0.5$. Then $\min\{|L|: L\in[I]^{<\infty},\ \delta$-dense$\}=2$. Let $K=\{0.5,0.95\}$. Clearly $d(K,I)=0.5$ and $K$ is not stretched.
\end{ex}

%------------------------------------------------------------------------------------------------------------------------------Unordered average---------------------------
\section{Unordered average}

\begin{df}A mean ${\cal K}$ of any number of variables is called \textbf{regular} if it satisfies the following three conditions: 
\begin{enumerate}
\item[(a)] ${\cal K}$ is permutation invariant i.e $\sigma:\{1,\dots,n\}\to\{1,\dots,n\}$ being a bijection implies that ${\cal K}(a_1, \dots, a_n)={\cal K}\big(a_{\sigma(1)}, \dots, a_{\sigma(n)}\big)$,

\item[(b)] $a_n\to \alpha\in\overline{\mathbb{R}}$ then ${\cal K}(a_1, \dots, a_n)\to\alpha$,

\item[(c)] if $J$ is an interval and ${\cal K}(a_1, \dots, a_n)\in J$ and $a_{n+1}\in J$ then 
\newline ${\cal K}(a_1, \dots, a_{n+1})\in J$ too.\qed
\end{enumerate}

\end{df}

Obviously the arithmetic mean \A\ (or almost every well known mean) satisfies these properties. Note that condition (c) does not imply condition (b) -- as it may seem at first glance (e.g. observe ${\cal K}=\min$).

\begin{df}\label{unorderedavg}Let $f:I\to\mathbb{R}$. Then we write ${\cal K}(\{f(i):i\in I\})=c\in\mathbb{R}$ if $\forall\epsilon>0$ there exists a finite set $I_0\subset I$ such that for every finite set $I_1$, $I_0\subset I_1\subset I$ implies that $\left|{\cal K}(\{f(i):i\in I_1\})-c\right|<\epsilon$. Then $c$ is called the \textbf{unordered average} of $f$ regarding ${\cal K}$.
\end{df}

The following observations show that in order to the unordered average exists, $f$ has to have a fairly simple structure.

\begin{prp}\label{punordavg1}If $f:I\to\mathbb{R}$ is given, the mean ${\cal K}$ is regular and the unordered average of $f$ regarding ${\cal K}$ exists then there is $c\in\mathbb{R}$ such that $f(i)=c$ except for countably many elements in $I$.
\end{prp}
\begin{proof}Suppose that there are two sequences $(a_n),\ (b_n)$ in $I$ such that $a_n\ne a_m,\ b_n\ne b_m\ (n\ne m),\ f(a_n)\to a,\ f(b_n)\to b$ and $a<b$. Suppose first that $a,b\in\mathbb{R}$. We can assume that $a_n\ne b_m$. Now if $I_0\in[I]^{<\infty}$ then when we add the elements of $I_0$ to $(a_n)$ we get $(a'_n)$, but $f(a'_n)\to a$ hence ${\cal K}(f(a'_1), \dots, f(a'_n))\to a$. Similarly when we add the elements of $I_0$ to $(b_n)$ we get $(b'_n)$, but $f(b'_n)\to b$ hence ${\cal K}(f(b'_1), \dots, f(b'_n))\to b$. Now if $\epsilon=\frac{b-a}{2}$ then we can create $I_1$ from the elements of $(a'_n)$ such that $I_0\subset I_1\in[I]^{<\infty}$ and ${\cal K}(\{f(i):i\in I_1\})<a+\epsilon$. Similarly we can create $I_2$ from the elements of $(b'_n)$ such that $I_0\subset I_2\in[I]^{<\infty}$ and ${\cal K}(\{f(i):i\in I_2\})>b-\epsilon$ which is a contradiction. 

If any of $a,b$ (or both) is infinite then a similar argument leads to a contradiction. 
\end{proof}

\begin{prp}\label{punordavg2}If $f:I\to\mathbb{R}$ is given, the mean ${\cal K}$ is regular and the unordered average of $f$ regarding ${\cal K}$ exists then $\{f(i):i\in I\}$ can have only one accumulation point that is the unordered average itself.
\end{prp}
\begin{proof}Actually in the proof of \ref{punordavg1} we showed the first part.

If one supposes that the accumulation point is not the unordered average then a similar argument that is in the proof of \ref{punordavg1} leads to a contradiction.
\end{proof}

\begin{thm}\label{tunordavg}Let $f:I\to\mathbb{R}$ be given and let the mean ${\cal K}$ be regular. Then the unordered average of $f$ regarding ${\cal K}$ exists and equals to $c$ iff either $\{i:f(i)\ne c\}$ is finite or $(f(i))_{i:f(i)\ne c}$ is a convergent sequence to $c$ where $c$ is the accumulation point of $\{f(i):i\in I\}$.
\end{thm}
\begin{proof}The necessity follows from \ref{punordavg1} and \ref{punordavg2}.

To prove the sufficiency if $K=\{i:f(i)\ne c\}$ is finite and $I$ is infinite then take a sequence $(a_n)$ in $I$ such that put the elements of $K$ to the beginning of the sequence. Then $f(a_n)\to c$ implies that ${\cal K}(f(a_1), \dots, f(a_n))\to c$. For given $\epsilon>0$ there is $N\geq |K|$ such that $|{\cal K}(f(a_1), \dots, f(a_N))-c|<\epsilon$. Let $I_0=\{a_1, \dots, a_N\}$. 

If $(f(i))_{i:f(i)\ne c}$ is a convergent sequence to $c$ and $\epsilon>0$ then choose $N$ such that $n\geq N$ implies that $|{\cal K}(f(a_1), \dots, f(a_n))-c|<\epsilon$ and $|f(a_n)-c|<\epsilon$. Now using condition (c) it is easy to check that $I_0=\{a_1, \dots, a_N\}$ satisfies the criterion in the definition of unordered average.
\end{proof}

\begin{cor}If there is $c\in\mathbb{R}$ such that the unordered sum exists for $f-c$ then the unordered average exists for $f$. \qed
\end{cor}

\begin{rem}If $I$ is a fully ordered set then we can define ${\cal K}(\{f(i):i\in I_1\})$ using the inheritance of the order from the finite set $I_1$.
Then one can leave condition (a) out and then develop the notion and properties of such unordered averages.
\end{rem}

Let us remark that neither ${\cal{M}}^{iso}(H)$ is an unordered average on the set of isolated points, nor ${\cal{M}}^{eds}(H)$ is an unordered average of $H$.

%----------------------------------------------------------------------------------------------------------------------The definitions of extensions--------------
\section{The definition of extensions}

\begin{df}\label{dext}Let $\langle I,d\rangle$ be a totally bounded pseudo-metric space. Let $s:[I]^{<\infty}\to\mathbb{R}$. The \textbf{extension of $\boldsymbol s$ onto $\boldsymbol I$} equals to $A$ if $\forall\epsilon>0\ \exists\delta>0$ such that $K\in[I]^{<\infty}$ and $d_H(K,I)<\delta$ implies that $|s(K)-A|<\epsilon$. In this case we use the notation ext$(s,I)=A$.
\end{df}

\begin{df}\label{dextstr}Let $\langle I,d\rangle$ be a totally bounded pseudo-metric space. Let $s:[I]^{<\infty}_{str}\to\mathbb{R}$. The \textbf{stretched extension of $\boldsymbol s$ onto $\boldsymbol I$} equals to $A$ if $\forall\epsilon>0\ \exists\delta>0$ such that $K\in[I]^{<\infty}_{str}$ and $d_H(K,I)<\delta$ implies that $|s(K)-A|<\epsilon$. In this case we use the notation ext$_{str}(s,I)=A$.
\end{df}

Clearly ext$(s,I)$ and ext$_{str}(s,I)$ are a kind of continuous extensions of $s$ to $I$ regarding the Hausdorff pseudo-metric.

\begin{rem}If in the previous definitions we replace $\mathbb{R}$ with any topological commutative group then we end up with a straightforward generalization of these notions. We are not going to investigate such a generalization further. 
\end{rem}

\begin{rem}\label{rgmtop}Another generalization can be the following. There is no metric/topology on $I$, instead we have a topology on the set $[I]^{<\infty}\cup\{\omega\}$ where $\omega\notin[I]^{<\infty}$ i.e. the points of this space are the finite subsets of $I$ and $\omega$. If $s:[I]^{<\infty}\to\mathbb{R}$ is given and $s$ can be extended to $\omega$ such that $s$ is continuous at $\omega$ then we say that ext$(s,I)$ exists and equals to the value on $\omega$.
 We are not going to investigate such a generalization further either. 
\end{rem}

\begin{df}\label{dextinfty}Let $\langle I,d\rangle$ be a totally bounded pseudo-metric space. Let $s:[I]^{<\infty}\to\mathbb{R}\ (s:[I]^{<\infty}_{str}\to\mathbb{R})$. The (stretched) extension of $s$ onto $I$ equals to $\infty$ if $\forall n\in\mathbb{N}\ \exists\delta>0$ such that $K\in[I]^{<\infty}$ (stretched) and $d_H(K,I)<\delta$ implies that $s(K)>n$.
\end{df}

One can similarly define when ext$(s,I)=-\infty$ (ext$_{str}(s,I)=-\infty$).

\begin{thm}\label{tequivdef}Let $\langle I,d\rangle$ be a totally bounded pseudo-metric space and let $d_H$ denote the Hausdorff pseudo-metric on the subsets of $I$. Then the following statements are equivalent:
\begin{enumerate}

\item[(1)] ext$(s,I)=A\in\mathbb{R}$

\item[(2)] $\forall\epsilon>0\ \exists\delta>0$ such that if $K\in[I]^{<\infty}$, $\delta$-sdense set in $I$, then $|s(K)-A|<\epsilon$

\item[(3)] $\forall\epsilon>0\ \exists\delta>0$ such that if $K,L\in[I]^{<\infty}$ and $d_H(K,I)<\delta,\ d_H(L,I)<\delta$, then $|s(K)-s(L)|<\epsilon$

\item[(4)] $\forall\epsilon>0\ \exists\delta>0$ such that if $K,L\in[I]^{<\infty}$, $\delta$-sdense sets in $I$, then $|s(K)-s(L)|<\epsilon$

\item[(5)] If $(K_n)$ is a sequence of subsets of $I$ such that $K_n\in[I]^{<\infty}$ and $d_H(K_n,I)\to 0$ then $\big(s(K_n)\big)$ is convergent.

\end{enumerate}
\end{thm}
\begin{proof}The equivalence of (1) and (2) and similarly (3) and (4) are straightforward from the already mentioned fact that $K$ is $\delta$-sdense iff $d_H(K,I)<\delta$.

Let us show (2)$\Rightarrow$(4). Choose $\delta$ for $\frac{\epsilon}{2}$ according to the definition. If $K,L\subset I$ are finite $\delta$-sdense sets then  $|s(K)-s(L)|\leq|s(K)-A|+|A-s(L)|<\frac{\epsilon}{2}+\frac{\epsilon}{2}=\epsilon$.

To prove (4)$\Rightarrow$(2) choose a finite $\frac{1}{n}$-sdense set $D_n$ for each $n\in\mathbb{N}$. By (4) $\big(s(D_n)\big)$ is a Cauchy sequence hence converges to a number, say $A$. For $\frac{\epsilon}{2}$ choose $\delta>0$ according to (4) and choose $n$ such that $|s(D_n)-A|<\frac{\epsilon}{2}$ and $\frac{1}{n}<\delta$. Now if $K\subset I$ is finite, $\delta$-sdense then $|s(K)-A|\leq|s(K)-s(D_n)|+|s(D_n)-A|<\frac{\epsilon}{2}+\frac{\epsilon}{2}=\epsilon$.

(1)$\Rightarrow$(5): Obviously $s(K_n)\to A$.

(5)$\Rightarrow$(1): Merging two such sequences shows that all such sequences converge to the same limit say $A$. Suppose that there is $\epsilon$ such that $\forall n\in\mathbb{N}$ there is $K_n\in[I]^{<\infty}$ such that $d_H(K_n,I)<\frac{1}{n},\ s(K_n)\notin(A-\epsilon,A+\epsilon)$. Then $(s(K_n))$ has a limit point different from $A$ -- a contradiction. 
\end{proof}

We can formulate a similar statement for stretched extension: juts add the adjective "stretched" in (2),(3),(4),(5). Also we can formulate statement version for $A=\pm\infty$.

We are going to investigate properties of ext$(s,I)$ and ext$_{str}(s,I)$ in details in section \ref{sprop}.

%------------------------------------------------------------------------------------------------------------------------------ Verifications  ---------------------------
\subsection{Verifications}

Now we are going to establish that the new notions cover most examples mentioned in the introduction. Arc length and integral on measure spaces will be handled later after we will have built some useful tools (see \ref{carclendefeq} and \ref{tintmeassp} respectively).

In all cases first we have to define a totally bounded pseudo-metric space $\langle I,d\rangle$ and a function $s:[I]^{<\infty}\to\mathbb{R}$ (or $s:[I]^{<\infty}_{str}\to\mathbb{R}$). Then we have to show that ext$(s,I)$ (or ext$_{str}(s,I)$) is equivalent to the old notion.

\subsubsection{Absolute convergent series}\label{ssabsconvser}

\begin{prp}\label{pserac}The following statements are equivalent:
\begin{enumerate}
\item[(1)] $\sum\limits_{i=1}^{\infty} a_i$ absolute convergent 

\item[(2)] $\forall\epsilon>0\ \exists N$ such that $\{1,2,\dots, N\}\subset A,B\subset\mathbb{N}$ being finite sets implies that $\Big|\sum\limits_{i\in A}a_i-\sum\limits_{i\in B}a_i\Big|<\epsilon$

\item[(3)] There is $c\in\mathbb{R}$ such that $\forall\epsilon>0\ \exists N$ such that $\{1,2,\dots, N\}\subset A\subset\mathbb{N}$ being a finite set implies that $\Big|\sum\limits_{i\in A}a_i-c\Big|<\epsilon$.
\end{enumerate}
\end{prp}
\begin{proof}(1)$\Rightarrow$(2), (1)$\Rightarrow$(3), (3)$\Rightarrow$(1) are straightforward.

We show (2)$\Rightarrow$(3). For all $n\in\mathbb{N}$ choose $D_n=\{1,2,\dots, N_n\}$ according to (2) and $\epsilon=\frac{1}{n}$ such that $n<m$ implies that $N_n<N_m$. Let $s_n=\sum\limits_{i\in D_n}a_i$. Then $\left(s_n\right)$ is a Cauchy sequence hence converges to say $c$. Let $\epsilon>0$. Choose $M$ such that $n>M$ implies that $|s_n-c|<\frac{\epsilon}{2}$. Choose $N_P>N_M$ such that $\{1,2,\dots, N_P\}\subset A,B\subset\mathbb{N}$ being finite sets implies that $\Big|\sum\limits_{i\in A}a_i-\sum\limits_{i\in B}a_i\Big|<\frac{\epsilon}{2}$. Then if $\{1,2,\dots, N_P\}\subset A$ then 
\[\left|\sum\limits_{i\in A}a_i-c\right|<\left|\sum\limits_{i\in A}a_i-\sum\limits_{i\in D_P}a_i\right|+\left|\sum\limits_{i\in D_P}a_i-c\right|<\frac{\epsilon}{2}+\frac{\epsilon}{2}=\epsilon.\qedhere\]
\end{proof}

If a series $\sum\limits_{i=1}^{\infty} a_i$ is given then let 

$I=\left\{\frac{1}{n}:n\in\mathbb{N}\right\}$ equipped with the usual metric,

$f:I\to\{a_n:n\in\mathbb{N}\},\ f\left(\frac{1}{n}\right)=a_n,$

$s:[I]^{<\infty}\to\mathbb{R},\ s(H)=\sum\limits_{h\in H}f(h)=\sum\limits_{h\in H}a_\frac{1}{h}\ \ \ (H\in[I]^{<\infty})$.

\begin{thm}Let an absolute convergent series $\sum\limits_{i=1}^{\infty} a_i$ be given. Then $\sum\limits_{i=1}^{\infty} a_i=c\in{\mathbb{R}}$ iff ext$(s,I)=c$ where $I,s$ are the same defined above.
\end{thm}
\begin{proof}Let $c\in\mathbb{R}$. First let $\sum\limits_{i=1}^{\infty} a_i=c$ and $\epsilon>0$. Then by \ref{pserac}(3) there is $N\in\mathbb{N}$ such that $\{1,2,\dots, N\}\subset A\subset\mathbb{N}$ being a finite set implies that $|\sum\limits_{i\in A}a_i-c|<\epsilon$. Now if $\delta=\frac{1}{N(N-1)}$ and $H\subset I$ is finite and $\delta$-sdense then $\left\{1,\frac{1}{2},\dots, \frac{1}{N}\right\}\subset H$ therefore $|s(H)-A|<\epsilon$  since $s(H)=\sum\limits_{i\in H}a_\frac{1}{i}$.

To prove the converse if ext$(s,I)=A$ and $\epsilon>0$ is given, choose $\delta>0$ according to the definition. Then let $N=\frac{1}{\delta}$. Now if $n>N$ then $\left\{\frac{1}{m}: m\leq n\right\}\subset L$ is $\delta$-sdense. Hence $|s(L)-A|<\epsilon$. But $s(L)=\sum\limits_{i=1}^{n} a_i$ which completes the proof.
\end{proof}

Similar statements can be formulated when $\sum a_i=+\infty$ and $\sum |a_i|^-<\infty$. We just formulate then and omit the proofs as they are almost the same. 

\begin{prp}\label{pserac2}The following statements are equivalent:
\begin{enumerate}
\item[(1)] $\sum\limits_{i=1}^{\infty} a_i=+\infty$ and $\sum\limits_{i=1}^{\infty} |a_i|^-<\infty$ 

\item[(2)] $\forall l\in\mathbb{R}\ \exists N$ such that $\{1,2,\dots, N\}\subset A\subset\mathbb{N}$ being a finite set implies that $l<\sum\limits_{i\in A}a_i$.\qed
\end{enumerate}
\end{prp}

\begin{thm}Let a series $\sum\limits_{i=1}^{\infty} a_i$ be given. Then the following statements are equivalent:
\begin{enumerate}
\item[(1)] $\sum\limits_{i=1}^{\infty} a_i=+\infty$ and $\sum\limits_{i=1}^{\infty} |a_i|^-<\infty$ 

\item[(2)] ext$(s,I)=+\infty$ where $s,I$ are the same defined above.\qed
\end{enumerate}
\end{thm}

\subsubsection{Unordered sum}

Let $I$ be an infinite set and $a:I\to\mathbb{R}$. Then let us define a pseudo-metric on $I$ in the following way: $d(x,y)=|a(x)-a(y)|\ (x,y\in I)$. Let $s:[I]^{<\infty}\to\mathbb{R},\ s(H)=\sum\limits_{h\in H}a(h)\ (H\in \text{Dom }s)$.

\begin{thm}Let $a:I\to\mathbb{R}$. Then the unordered sum exists and equals to $A$ iff ext$(s,I)=A$  where $I,s$ are the same defined above.
\end{thm}
\begin{proof}Starting with necessity it is known that $a(I)$ has only one accumulation point that is 0. It follows that $a^{-1}(\mathbb{R}-\{0\})$ is countable and the set of its accumulation points is $a^{-1}(\{0\})$. Let $\epsilon>0$. Then there is $I_0\subset I$ finite such that $I_0\subset J\subset I$ finite implies that $|s(J)-A|<\epsilon$. We can assume that $I_0\cap a^{-1}(\{0\})=\emptyset$. We can also assume that $I_0$ is closed since simply replace $I_0$ with $cl(I_0)$ (that is finite too). Set $\delta=d(I_0,I-I_0)>0$. Clearly $K\subset I$ being $\frac{\delta}{2}$-sdense implies that $I_0\subset K$ that gives that $|s(K)-A|<\epsilon$. Which yields that ext$(s,I)=A$.

To show the sufficiency let $\epsilon>0$. Then there is $\delta>0$ such that $K\subset I$ being finite and $\delta$-sdense implies $|s(K)-A|<\epsilon$. Choose an arbitrary  $I_0\subset I$ finite and $\delta$-sdense set. Then obviously $I_0\subset J\subset I$ being finite implies that $J$ is $\delta$-sdense too hence that $|s(J)-A|<\epsilon$ which completes the proof.
\end{proof}

\begin{rem}According to \ref{rgmtop} we can set $s$ as above and take the following topology on $[I]^{<\infty}\cup\{\omega\}$: let $[I]^{<\infty}$ be a discrete subspace and for $\omega$ take the cofinite neighborhood filter. Then it can be easily seen that ext$(s,I)$ exists iff the unordered sum exists and they are equal.
\end{rem}

\subsubsection{Series in general}

If a series $\sum\limits_{i=1}^{\infty} a_i$ is given then let 

$I=\left\{\frac{1}{2^n}:n\in\mathbb{N}\right\}$ equipped with the usual metric,

$f:I\to\{a_n:n\in\mathbb{N}\},\ f\left(\frac{1}{2^n}\right)=a_n,$

$s:[I]^{<\infty}\to\mathbb{R},\ s(H)=\sum\limits_{h\in H}f(h)=\sum\limits_{h\in H}a_{-\log_2 h}\ (H\in \text{Dom }s)$.

\begin{prp}\label{psertn}Let $I=\left\{\frac{1}{2^n}:n\in\mathbb{N}\right\}$. Let $K\subset I$ be finite, stretched, $|K|\geq 2$. Then there is $k_0,k_1\in K$ such that $K=\left\{i\in I:i \geq k_1\right\}\cup\{k_0\}$. Actually $k_0=\min K, k_1=\min K-\{k_0\}$.
\end{prp}
\begin{proof}Let $k_0=\min K, k_1=\min K-\{k_0\}$. Suppose that there is $i\in I$ such that $i>k_1$ and $i\notin K$. Let 
\[m=\min\{i\in I:i>k_1,\ i\notin K\},\ k_2=\max\{k\in K:k<m\}.\]
Clearly $m\notin K,\ k_2\in K,\ k_0<k_1\leq k_2<m$ and $(k_2,m)\cap I=\emptyset$. This yields that $d_H(K,I)\geq k_2$. But $d(k_0,k_1)<k_1\leq k_2\leq d_H(K,I)$ which contradicts to stretchedness.
\end{proof}

\begin{prp}\label{psertn2}Let $I=\left\{\frac{1}{2^n}:n\in\mathbb{N}\right\}$. Let $K\subset I$ be finite, stretched, $d_H(K,I)<\frac{1}{2^N}$. Then 
\[\left\{i\in I:i \geq \frac{1}{2^{N-1}}\right\}\subset K.\]
\end{prp}
\begin{proof}By \ref{psertn} $K=\left\{i\in I:i \geq k_1\right\}\cup\{k_0\}$ where $k_0=\min K, k_1=\min K-\{k_0\}$. If we supposed that $\frac{1}{2^{N-1}}\notin K$ then it would give that $d_H(K,I)\geq\frac{1}{2^N}$ which is a contradiction. Clearly $k_0\leq\frac{1}{2^N}$ which yields that $k_1\leq\frac{1}{2^{N-1}}$.
\end{proof}

\begin{thm}Let a series $\sum\limits_{i=1}^{\infty} a_i$ be given. Then $\sum\limits_{i=1}^{\infty} a_i=A$ iff ext$_{str}(s,I)=A$ where $I,s$ are the same defined above.
\end{thm}
\begin{proof}Let $A\in\mathbb{R}$ first.

Let $\sum\limits_{i=1}^{\infty} a_i=A$ and $\epsilon>0$. Then there is $N$ such that $n\geq N$ implies that  $\left|\sum\limits_{i=1}^{n} a_i-A\right|<\frac{\epsilon}{2}$ and $|a_n|<\frac{\epsilon}{2}$. Now let $\delta=\frac{1}{2^{N+2}}$. Because of \ref{psertn2} if $K\subset I$ is finite, stretched and $\delta$-sdense then $K=\left\{i\in I:i \geq k_1\right\}\cup\{k_0\}$ where $k_0=\min K\leq \frac{1}{2^{N+1}},\ k_1=\min K-\{k_0\}\leq\frac{1}{2^{N}}$. Hence 
\[s(K)=a_{-\log_2k_0}+\sum\limits_{i=1}^{-\log_2k_1} a_i.\] 
Therefore 
\[|s(K)-A|\leq|a_{-\log_2k_0}|+\left|\sum\limits_{i=1}^{-\log_2k_1} a_i-A\right|<\frac{\epsilon}{2}+\frac{\epsilon}{2}=\epsilon\]
which yields that ext$_{str}(s,I)=A$.

\smallskip

Now if ext$_{str}(s,I)=A$ and $\epsilon>0$ is given then choose $\delta>0$ according to the definition. Take a $K\subset I$ finite set that is stretched and $\delta$-sdense. Then let $N=-\log_2\min K.$ Now if $n>N$ then $L=\left\{\frac{1}{2^m}: m\leq n\right\}$ is stretched and $\delta$-sdense. Hence $|s(L)-A|<\epsilon$. But $s(L)=\sum\limits_{i=1}^{n} a_i$ which completes the proof.

\smallskip

If $A=+\infty\ (-\infty)$ then $\sum\limits_{i=1}^{\infty} a_i=+\infty\ (-\infty)$ gives that $(a_i)$ bounded from below (above). Using that fact a similar argument can be applied.
\end{proof}

\begin{thm}Let $\sum\limits_{i=1}^{\infty} a_i$ be convergent but non absolute convergent series, converging to $A$. Then there are no totally bounded pseudo-metric space $I,\ f:I\to\{a_n:n\in\mathbb{N}\}$ surjective mapping, and $s:[I]^{<\infty}\to\mathbb{R},\ s(H)=\sum\limits_{h\in H}f(h)\ (H\in \text{Dom }s)$ such that ext$(s,I)=A$.
\end{thm}
\begin{proof}Suppose the contrary and for $\epsilon=1$ choose $\delta$ such that $K\in[I]^{<\infty}$ $\delta$-sdense implies that $|s(K)-A|<1$. Choose a $\delta$-sdense $K_0\in[I]^{<\infty}$. If $\sum\limits_{i=1}^{\infty} a_i$ is non absolute convergent then there is a subsequence $(a_{n_j})$ such that $\sum\limits_{j=1}^{\infty} a_{n_j}=+\infty$. For each $j\in\mathbb{N}$ choose $i_j\in I$ such that $f(i_j)=a_{n_j}$. Now choose elements from $\{i_j:j\in\mathbb{N}\}-K_0$ such that for the set $K_1$ of those elements, $s(K_1)>2$ holds. Then $K_0\cup K_1$ is still $\delta$-sdense but $|s(K_0\cup K_1)-A|<1$ does not hold which is a contradiction.
\end{proof}

\subsubsection{Riemann integral}\label{ssriemann}

If $f:[a,b]\to\mathbb{R}$ is given then let $I=[a,b]$ equipped with the usual metric and $s:[I]^{<\infty}\to\mathbb{R}$ defined in the following way. Let $H\subset I$ be finite. Let $\tilde H=H\cup\{a,b\}$. Let $\tilde H=\{a_0,\dots,a_n\}$ where $a=a_0<\dots <a_n=b$.

Now let \[s(H)=f(a_1)(a_2-a_0)+f(a_3)(a_4-a_2)+\dots+
\begin{cases}
f(a_{n-1})(a_{n}-a_{n-2}) &\text{if }n\text{ is even}\\  
f(b)(a_{n}-a_{n-1}) & \text{if }n\text{ is odd.}
\end{cases}
\]

It is easy to see that if $H$ is $\delta$-sdense then so is $\tilde H$. It is also clear that if $H$ is $\delta$-sdense then $|a_i-a_{i-1}|<2\delta\ (1\leq i\leq n)$.

\begin{thm}Let a bounded $f:[a,b]\to\mathbb{R}$ be given. Then $\int\limits_a^bf=A$ iff ext$(s,I)=A$ where $s,I$ are the same defined above.
\end{thm}
\begin{proof}To see the necessity let $\epsilon>0$. Then there is $\delta>0$ such that the Riemann sum is in $\epsilon$ distance within $A$ for any subintervals of length less than $\delta$ and for any chosen points in between. Therefore if we take a $\frac{\delta}{4}$-sdense set $H$ on $I$ then $\tilde H=\{a_0,\dots,a_n\}$ is $\frac{\delta}{4}$-sdense too. Hence if $n$ is even and we take $[a_0,a_2],[a_4,a_2],\dots [a_n,a_{n-2}]$ as subintervals and $a_1,a_3,\dots a_{n-1}$ as chosen points then $s(H)$ will be the associated Riemann sum. Similarly if $n$ is odd and we take $[a_0,a_2],[a_4,a_2],\dots [a_n,a_{n-1}]$ as subintervals and $a_1,a_3,\dots a_{n-2},b$ as chosen points then $s(H)$ will be the associated Riemann sum. In both cases the length of subintervals is less than $\delta$. That yields that $|s(H)-A|<\epsilon$ which means that ext$(s,I)=A$. 

\smallskip

To show the sufficiency let again $\epsilon>0$ and $|f|<M$. Then there is $\delta>0$ such that $2\delta M<\frac{\epsilon}{3}$ and $|s(H)-A|<\frac{\epsilon}{3}$ whenever $H$ is $\delta$-sdense and finite. Let us take points $a=a_0<a_1<\dots<a_n=b$ such that $a_i-a_{i-1}<\frac{\delta}{2}\ (1\leq i\leq n)$. Then choose points $\xi_i\in[a_{i-1},a_i]\ (1\leq i\leq n)$. 

It may happen for two subintervals next to each other that the chosen point is the same: for $[a_{i-2},a_{i-1}]$ and $[a_{i-1},a_i]$  $\xi_{i-1}=\xi_i=a_{i-1}$ holds. Then drop $a_{i-1}$ from the list, i.e join the two intervals into $[a_{i-2},a_i]$ and keep $\xi_i$ and drop the other. Now apply this operation for all such pairs of subintervals. With this operation we end up with a subinterval structure containing less subintervals and having the same Riemann sum but two consecutive intervals do not have the same chosen point.

Now we are going to consider the Riemann sum on all subintervals except the first and last subintervals and we are going to modify it slightly.
Count those subintervals where $\xi_i=a_{i-1}$ or $a_i$  (and $i\ne 1,\ i\ne n$) holds. Say there are $m$ such subintervals. If for such subinterval $\xi_i=a_{i}$ holds then $\xi_{i+1}>\xi_i$ (by the result of the previous operation). Then it is readily seen that it is possible to move upwards $a_i$ slightly such that $\xi_i\in(a_{i-1},a_i)$ and $\xi_{i+1}\in(a_{i},a_{i+1})$ holds and the Riemann sum  changes with less than $\frac{\epsilon}{3m}$. The similar other case when $\xi_i=a_{i-1}$ can be handled similarly.
If we apply this second operation for all $m$ subintervals then we end up with a new subinterval structure where the chosen points are never equal to any of the end points of the subintervals (except for the first and last subintervals) and the new Riemann sum excluding the first and last subintervals differs from the old one by less than $\frac{\epsilon}{3}$.  

As a third operation we modify the the first and last subintervals midpoints if needed. If any of the chosen points of the first and last subintervals is equal to any of the endpoints of the subinterval then replace the midpoint with a new one that is not equal to any of the endpoints of the subinterval.
Note that the contribution of the sums on the first and last subintervals is less than $\frac{\epsilon}{3}$ altogether. Hence this new Riemann sum differs from the original one by less than $\frac{2\epsilon}{3}$.

Let $H$ contain all end points of the subintervals and the chosen points as well. Then $s(H)$ equals to the Riemann sum. Clearly $H$ is $\delta$-sdense hence $|s(H)-A|<\frac{\epsilon}{3}$. But the modified Riemann sum differs from the original one by less than $\frac{2\epsilon}{3}$ therefore the original Riemann sum differs from $A$ by less than $\epsilon$.
\end{proof}

\subsubsection{Unordered average}

Let ${\cal K}$ be a regular mean, $I$ be an infinite set and $a:I\to\mathbb{R}$. Then let us define a pseudo-metric on $I$ in the following way: $d(x,y)=|a(x)-a(y)|\ (x,y\in I)$. Let 
\[s:[I]^{<\infty}\to\mathbb{R},\ s(K)={\cal K}(\{a(i):i\in K\})\ \ \ (K\in \text{Dom }s).\]

\begin{thm}Let $a:I\to\mathbb{R}$ and $\{a(i):i\in I\}$ is infinite. Then the unordered average exists and equals to $c$ iff ext$(s,I)=c$  where $I,s$ are the same defined above.
\end{thm}
\begin{proof}Starting with necessity we know from \ref{tunordavg} that either $\{i:a(i)\ne c\}$ is finite or $(a(i))_{i:a(i)\ne c}$ is a convergent sequence to $c$.

Let $\{i\in I:a(i)\ne c\}=\{i_n\in I:n\in\mathbb{N}\}$ such that $(a(i_n))$ is a convergent sequence to $c$ and let $\epsilon>0$. Then by condition (b) there is $N$ such that $n\geq N$ implies that $|{\cal K}(a(i_1),\dots,a(i_n))-c|<\epsilon$ and $|a(i_n)-c|<\epsilon$. Let $I_0=\{i_1,\dots,i_N\}$. We can also assume that $I_0$ is closed since simply replace $I_0$ with $cl(I_0)$ (that is finite too). Set $\delta=d(I_0,I-I_0)>0$. Clearly $J\in[I]^{<\infty}$ being $\frac{\delta}{3}$-sdense implies that $I_0\subset J$. If $j\in J-I_0$ then $|a(j)-c|<\epsilon$ which gives that $|{\cal K}(\{a(j):j\in J\})-c|<\epsilon$ by condition (c). That is $|s(J)-c|<\epsilon$ which yields that ext$(s,I)=c$.

To show the sufficiency let $\epsilon>0$. Then there is $\delta>0$ such that $K\subset I$ being finite and $\delta$-sdense implies $|s(K)-c|<\epsilon$. Choose an arbitrary  $I_0\subset I$ finite and $\delta$-sdense set. Then obviously $I_0\subset J\subset I$ being finite implies that $J$ is $\delta$-sdense too hence that $|s(J)-c|<\epsilon$ which completes the proof.
\end{proof}

\begin{ex}If $\{a(i):i\in I\}$ is finite then a straightforward example shows that the necessity part of the theorem is false. Let $I=\big\{\frac{1}{n}:n\in\mathbb{N}\big\},\ a(1)=1, a(i)=0$ if $i\in I,i\ne 1$. Clearly the unordered average exists and equals to $0$. If $K=\big\{1,\frac{1}{2}\big\}$ then $K$ is $\delta$-sdense for all $\delta\leq 1$ and $s(K)=\frac{1}{2}$. Similarly if $L=\big\{1,\frac{1}{2},\frac{1}{3}\big\}$ then $L$ is $\delta$-sdense for all $\delta\leq 1$ and $s(L)=\frac{1}{3}$. Which show that ext$(s,I)$ does not exists.
\end{ex}

\subsubsection{$\boldsymbol{{\cal{M}}^{iso}}$}

\begin{lem}\label{liso1}Let a bounded sequence $(a_i)$ be given such that $\forall i\ |a_i|<M$. Then $|\A(a_1,\dots,a_{n})-\A(a_1,\dots,a_{n+1})|<\frac{2M}{n+1}$.
\end{lem}
\begin{proof}\[|\A(a_1,\dots,a_{n})-\A(a_1,\dots,a_{n+1})|=\left|\frac{\sum\limits_{i=1}^na_i}{n}-\frac{\sum\limits_{i=1}^{n}a_i+a_{n+1}}{n+1}\right|=\]
\[\left|\frac{n\sum\limits_{i=1}^na_i+\sum\limits_{i=1}^na_i-n\sum\limits_{i=1}^na_i-na_{n+1}}{n(n+1)}\right|=\left|\frac{\frac{\sum\limits_{i=1}^na_i}{n}-a_{n+1}}{n+1}\right|\leq\frac{\left|\frac{\sum\limits_{i=1}^na_i}{n}\right|+|a_{n+1}|}{n+1}<\frac{2M}{n+1}.\qedhere\]
\end{proof}

\begin{lem}\label{liso2}Let a bounded sequence $(a_i)$ be given such that $\forall i\ |a_i|<M$. Let $n<m,\ A=\A(a_1,\dots,a_{n}),\ B=\A(a_1,\dots,a_{m})$. Then there is $\sigma:\{n+1,\dots,m\}\to\{n+1,\dots,m\}$ permutation such that 
\[\forall k\in\{n+1,\dots,m\}\ \ \A(a_1,\dots,a_{n},a_{\sigma(n+1)},\dots,a_{\sigma(k)})\in \left(A-\frac{2M}{n+1},B+\frac{2M}{n+1}\right)
\]
if $A\leq B$ (if $B<A$ then swap them).
\end{lem}
\begin{proof}We define $\sigma$ recursively. Because of $A\leq B$ there is $l\in\{n+1,\dots,m\}$ such that $a_l\geq A$. Let $\sigma(n+1)=l$. By \ref{liso1} it satisfies the condition.

Suppose we have already defined $\sigma$ for $n+1,\dots,k$ and 
\[A_k=\A(a_1,\dots,a_{n},a_{\sigma(n+1)},\dots,a_{\sigma(k)})\in \left(A-\frac{2M}{n+1},B+\frac{2M}{n+1}\right).\]
There are 3 cases:
\begin{enumerate}
\item If $A_k\leq A$ then there is $l\in\{n+1,\dots,m\}$ that is not chosen yet such that $a_l\geq A$. Let $\sigma(k+1)=l$.

\item If $A_k\in(A,B)$ then choose arbitrary $l\in\{n+1,\dots,m\}$ that is not chosen yet and let $\sigma(k+1)=l$.

\item If $A_k\geq B$ then there is $l\in\{n+1,\dots,m\}$ that is not chosen yet such that $a_l\leq B$. Let $\sigma(k+1)=l$.
\end{enumerate}
In each case \ref{liso1} guaranties that the condition will hold.
\end{proof}

Let $H\subset\mathbb{R}$ be bounded, infinite and $cl(H-H')=H$ and $h\in H$ implies that $|h|<M$. Now we are going to define $I$ and $s$.

\smallskip

Let $D=\{d(x,H'):x\in H-H'\}=\{t_1,t_2\dots\}$ such that $t_1>t_2>\dots$ Let $D_n=\{x\in H-H':d(x,H')=t_n\}\ (n\in\mathbb{N})$. Clearly $D_n$ is finite and let $|D_n|=e_n\ (n\in\mathbb{N})$ and $g_n=\sum\limits_{i=1}^{n}e_n$. 

\smallskip

Now we are going to put the elements of $H-H'$ into a sequence $(a_i)$ which is possible hence $H-H'$ is countably infinite. Arrange the elements of $D_1$ into $(a_i)$ arbitrarily. Then we go on by recursion. Suppose we are done till $n\in\mathbb{N}$ and all elements of $D_1,\dots,D_n$ have been added to $(a_i)$. Let $A_n=\A(D_1\cup\dots\cup D_n),\ A_{n+1}=\A(D_1\cup\dots\cup D_{n+1})$. Then add the elements of $D_{n+1}$ to $(a_i)$ according to \ref{liso2} i.e. if $g_n<j\leq g_{n+1}$ then $a_j\in\left(A_n-\frac{2M}{g_n+1},A_{n+1}+\frac{2M}{g_n+1}\right)$ if $A_n\leq A_{n+1}$ (if $A_{n+1}<A_n$ then swap them)

\smallskip

$I=\left\{\frac{1}{2^n}:n\in\mathbb{N}\right\}$ equipped with the usual metric,

$f:I\to\{a_n:n\in\mathbb{N}\},\ f\left(\frac{1}{2^n}\right)=a_n,$

Let $s:[I]^{<\infty}_{str}\to\mathbb{R},\ s(K)=\A(f(K))=\A(\{f(k):k\in K\})\ (K\in \text{Dom }s)$.

\begin{thm}Let $H\subset\mathbb{R}$ be bounded and $cl(H-H')=H$. Then ${\cal{M}}^{iso}(H)=A$ iff ext$_{str}(s,I)=A$ where $I,s$ are the same defined above.
\end{thm}
\begin{proof}First let ${\cal{M}}^{iso}(H)=A$ and $\epsilon>0$. Then there is $N\in\mathbb{N}$ such that $n\geq N$ implies that $\left|A_n-A\right|<\frac{\epsilon}{3}$ and $\frac{2M}{g_n}<\frac{\epsilon}{3}$. Let  $\delta=\frac{1}{2}\cdot\frac{1}{2^{g_N}}$.

Now if $K\subset I$ is finite, stretched and $\delta$-sdense then by \ref{psertn} there is $k_0,k_1\in K$ such that $K=\left\{i\in I:i \geq k_1\right\}\cup\{k_0\}$. Then there is $n\geq N$ such that $D_1\cup\dots\cup D_n\subset f(K-\{k_0\})\subset D_1\cup\dots\cup D_{n+1}$. By \ref{liso2} 
\[\A(f(K-\{k_0\}))\in\left(A_n-\frac{\epsilon}{3},A_{n+1}+\frac{\epsilon}{3}\right)\subset\left(A-\frac{2\epsilon}{3},A+\frac{2\epsilon}{3}\right)\]
if $A_n\leq A_{n+1}$ (if $A_{n+1}<A_n$ then swap them). By \ref{liso1} $|\A(f(K))-A|=|s(K)-A|<\epsilon$ which gives that ext$_{str}(s,I)=A$.

To prove the converse if ext$_{str}(s,I)=A$ and $\epsilon>0$ is given, choose $\delta>0$ according to the definition. Then choose $N$ such that  $\frac{1}{2^{g_N}}<\frac{1}{2}\delta$. Now if $\delta'<t_N$ then there is $n\geq N$ such that $H-S(H',\delta')=D_1\cup\dots\cup D_n$. But if $K=f^{-1}(D_1\cup\dots\cup D_n)$ then $K$ is $\delta$-sdense and stretched therefore $|s(K)-A|<\epsilon$ which yields that $|\A(f(K))-A|=|\A(H-S(H',\delta'))-A|<\epsilon$ which gives that ${\cal{M}}^{iso}(H)=A$.
\end{proof}

\subsubsection{$\boldsymbol{{\cal{M}}^{eds}}$}

Let $H\subset\mathbb{R}$ be bounded. 

Let $I=H$ and $s:[I]^{<\infty}\to\mathbb{R},\ s(K)=\A(K)\ \ (K\in \text{Dom }s)$.

\begin{thm}Let $H\subset\mathbb{R}$ be bounded. If ext$(s,I)=A$ then ${\cal{M}}^{eds}(H)=A$ where $I,s$ are the same defined above.
\end{thm}
\begin{proof}Let ext$(s,I)=A$ and $\epsilon>0$. Then there is $\delta>0$ such that $K\subset I$ being finite and $\delta$-sdense implies that $|s(K)-A|<\epsilon$. Let $a=\inf H,b=\sup H,\ N=\frac{b-a}{\delta}$. If $n>N, \delta'=\frac{b-a}{n}$ then clearly $0<\delta'<\delta$. Let $K\in[I]^{<\infty}$ such that
\[\forall i\in\mathbb{N}\cup\{0\}\ \ \Big|K\cap\big[a+i\delta',a+(i+1)\delta'\big)\Big|=\mathrm{sign}\bigg(\Big|H\cap\big[a+i\delta',a+(i+1)\delta'\big)\Big|\bigg)\]
where $\mathrm{sign}$ takes 1 for infinite cardinals too.
Obviously $K$ is $\delta$-sdense hence $|s(K)-A|<\epsilon$ which gives that ${\cal{M}}^{eds}(H)=A$.
\end{proof}

\begin{prp}Let $C$ be the ternary Cantor set. Then ext$_{str}(s,C)$ does not exist.
\end{prp}
\begin{proof}Consider the elements of $C$ in base 3 expansion i.e. 
\[C=\{x\in[0,1]: x={0.x_1x_2\dots}_{\textcircled{\tiny{3}}},\ x_i\in\{0,2\}\}.\]
For $n\in\mathbb{N}$ let 
\[K_n=\{x\in C: x={0.x_1x_2\dots}_{\textcircled{\tiny{3}}},\ x_{n}=x_{n+1}=\dots\}\]
and 
\[L_{1,n}=\{x\in C: x={0.0x_2x_3\dots}_{\textcircled{\tiny{3}}},\ x_{n}=x_{n+1}=\dots;\ x_{n-1}=2,x_n=0,\text{ or } x_{n-1}=0,x_n=2\}\]
\[L_{2,n}=\{x\in C: x={0.2x_2x_3\dots}_{\textcircled{\tiny{3}}},\ x_{n}=x_{n+1}=\dots\}\]
\[L_n=L_{1,n}\cup L_{2,n}.\]
Obviously $d_H(K_n,C)=d_H(L_n,C)=3^{-n}$ and both sets are stretched. Moreover $|K_n|=2^n,\ |L_{1,n}|=2^{n-2},\ |L_{2,n}|=2^{n-1}$. It is easy to check that $\A(K_n)=\frac{1}{2}$ while $L_{1,n}$ is symmetric to $\frac{1}{6}$ and $L_{2,n}$ is symmetric to $\frac{5}{6}$ which gives that 
\[\A(L_n)=\frac{2^{n-2}\cdot\frac{1}{6}+2^{n-1}\cdot\frac{5}{6}}{2^{n-2}+2^{n-1}}=\frac{\frac{1}{6}+2\cdot\frac{5}{6}}{1+2}=\frac{11}{18}\]
which shows that ext$_{str}(s,C)$ cannot exist.
\end{proof}

Therefore we have to find a different method in order to produce ${\cal{M}}^{eds}$ as an extension. For that we have to modify (reduce) the domain of $s$.

\smallskip

Let $I=H$ and ${\cal S}=\bigcup\limits_{n=1}^{\infty}{\cal S}_n$ where
\[{\cal S}_n=\Bigg\{K\in[I]^{<\infty}:\forall i\in\{0,\dots,n-1\}\]
\[\Big|K\cap\big[a+i\delta,a+(i+1)\delta\big)\Big|=\mathrm{sign}\bigg(\Big|H\cap\big[a+i\delta,a+(i+1)\delta\big)\Big|\bigg)\text{ where }\delta=\frac{b-a}{n}\Bigg\}\]
$\tilde{s}:{\cal S}\to\mathbb{R},\ s(K)=\A(K)\ \ (K\in {\cal S})$.

\begin{lem}\label{leds1}$\forall n\in\mathbb{N}\ \exists \delta_n>0$ such that $\forall K\in{\cal S}_n\ d_H(K,H)\geq\delta_n$.
\end{lem}
\begin{proof}Let $K\in{\cal S}_n$. As $H$ is infinite, there is $ i_0\in\{0,\dots,n-1\}$ such that $\exists x_1,x_2\in H\cap\big[a+i_0\frac{b-a}{n},a+(i_0+1)\frac{b-a}{n}\big);\ x_1<x_2;\ x_1,x_2\notin K$ and $x_1,x_2\ne a+i_0\frac{b-a}{n}$. Let $t_1=d\big(a+i_0\frac{b-a}{n},x_1\big),\ t_2=d(x_1,x_2),\ t_3=d\big(a+(i_0+1)\frac{b-a}{n},x_2\big)$. Then there are $y_1,y_2\in K$ such that $d(x_1,y_1)\leq d_H(K,H),d(x_2,y_2)\leq d_H(K,H)$. Then there are 3 cases:
\begin{enumerate}
\item If $y_1\notin\big[a+i_0\frac{b-a}{n},a+(i_0+1)\frac{b-a}{n}\big)$ then $d(x_1,y_1)>t_1$.
\item If $y_2\notin\big[a+i_0\frac{b-a}{n},a+(i_0+1)\frac{b-a}{n}\big)$ then $d(x_2,y_2)\geq t_3$.
\item If $y_1,y_2\in\big[a+i_0\frac{b-a}{n},a+(i_0+1)\frac{b-a}{n}\big)$ then $y_1=y_2$ and either $d(x_1,y_1)\geq \frac{t_2}{2}$ or $d(x_2,y_2)\geq \frac{t_2}{2}$.
\end{enumerate}
In either case $d_H(K,H)\geq\min\{t_1,t_3,\frac{t_2}{2}\}$.
\end{proof}

\begin{thm}Let $H\subset\mathbb{R}$ be bounded. If ${\cal{M}}^{eds}(H)=A$ then ext$(\tilde{s},I)=A$ where $I,\tilde{s}$ are the same defined above.
\end{thm}
\begin{proof}By \ref{leds1} for all $n\in\mathbb{N}$ there is $\delta_n>0$ such that $\forall K\in{\cal S}_n\ d_H(K,H)\geq\delta_n$.

For given $\epsilon>0$ choose $N$ according to the definition of  ${\cal{M}}^{eds}$ i.e. if $n>N,\ K\in {\cal S}_n$ then $|\tilde{s}(K)-A|<\epsilon$. Let $\delta=\min\{\delta_n:n\leq N\}$. Now if $K\in[I]^{<\infty},K\in{\cal S},\ d_H(K,H)<\delta$ then $\forall n\leq N\ K\notin {\cal S}_n$ hence $\exists m>N$ such that $K\in {\cal S}_m$. Therefore $|\tilde{s}(K)-A|<\epsilon$ which gives that ext$(\tilde{s},I)=A$.
\end{proof}

\subsubsection{Convergence of a sequence}

The convergence of a sequence can be fit into this framework however it does not seem to have any added value to the theory of convergence. However it is shown at least that the framework can handle that notion as well.

\smallskip

Let $(a_n)$ be a given sequence. 

Let $I=\left\{\frac{1}{n}:n\in\mathbb{N}\right\}$ and $s:[I]^{<\infty}\to\mathbb{R},\ s(K)=a_{\frac{1}{\min K}}\ \ \ (H\in \text{Dom }s)$.

\begin{thm}Let $(a_n)$ be a given sequence. Then $a_n\to A$ iff ext$(s,I)=A$ where $I,s$ are the same defined above.
\end{thm}
\begin{proof}If $a_n\to A\in\mathbb{R}, \epsilon>0$ then there is $N$ such that $n>N$ implies that $|a_n-A|<\epsilon$. If $K\in[I]^{<\infty}$ and $\frac{1}{N}$-sdense then $|s(K)-A|<\epsilon$.

To see the converse let $\delta$ be chosen to the definition of ext$(s,I)$. Let $N=\frac{1}{\delta}$. If $n>N$ then let $K=\{i\in I:i\geq\frac{1}{n}\}$. Clearly $K$ is $\delta$-sdense hence $|s(K)-A|<\epsilon$ which is equivalent to $|a_n-A|<\epsilon$.

\smallskip

$A=\pm\infty$ can be handled similarly.
\end{proof}

%----------------------------------------------------------------------------------------------------------------------Properties of ext$(s,I)$ and ext$_{str}(s,I)$s--------------
\section{Properties of ext$(s,I)$ and ext$_{str}(s,I)$}\label{sprop}

Throughout this section $\langle I,d\rangle$ will always denote a totally bounded pseudo-metric space and $s$ will denote a function $s:[I]^{<\infty}\to\mathbb{R}$.

%--------------------------------
\subsection{Basic properties}

\begin{prp}Both ext$(s,I)$ and ext$_{str}(s,I)$ are unique if exist.
\end{prp}
\begin{proof}Suppose the $A_1$ and $A_2$ both satisfies definition \ref{dext} (\ref{dextstr}). Then let $\epsilon=\frac{|A_1-A_2|}{2}$. Choose $\delta$ according to the definition and take $K\subset I$ that is finite and $\delta$-sdense (and stretched) to get a contradiction.
\end{proof}

\begin{prp} \hfill
\begin{enumerate}
\item [(a)] If ext$(s,I)$ exists then so does ext$_{srt}(s,I)$ and they are equal.

\item [(b)] Let ext$_{srt}(s,I)$ exist. If $\forall\epsilon>0\ \exists\delta>0$ such that $K\in[I]^{<\infty}_{str}\ d_H(K,I)<\delta,\ L\in[I]^{<\infty},\ d_H(L,I)<\delta$ implies that $|s(K)-s(L)|<\epsilon$ then ext$(s,I)$ exists and ext$(s,I)$=ext$_{srt}(s,I)$.
\end{enumerate}
\end{prp}
\begin{proof}(a) is straightforward.

(b): Apply \ref{tequivdef}(3) and then (a).
\end{proof}

\begin{prp}\label{punifcontofs}If $s$ is uniformly continuous on $[I]^{<\infty}$ regarding the Hausdorff pseudo-metric then ext$(s,I)$ exists and finite.  
\end{prp}
\begin{proof}Apply \ref{tequivdef}(3).
\end{proof}

\begin{ex}The opposite is not true. Let $I=\mathbb{Q}\cap(0,1)$ and let $s$ be defined in the way that for sets containing an only point $s$ takes 1, for sets containing more than one point $s$ takes 0. Obviously ext$(s,I)$ exists and equals to $0$ however $s$ is not even continuous.
\end{ex}

\begin{prp}\label{plinearity}Let $s_1:[I]^{<\infty}\to\mathbb{R},\ s_2:[I]^{<\infty}\to\mathbb{R}$ and both ext$(s_1,I)$ and ext$(s_2,I)$ exist and be finite. Let a continuous function $F:\mathbb{R}\times\mathbb{R}\to\mathbb{R}$ be given and let $s(K)=F(s_1(K),s_2(K))$ ($K\in[I]^{<\infty}$). Then ext$(s,I)$ exists and equals to $F($ext$(s,I_1),$ext$(s,I_2))$.  
\end{prp}
\begin{proof}Let ext$(s_1,I)=A_1$ and ext$(s_2,I)=A_2$. Let $\epsilon>0$ and choose $\epsilon_1$ according to the continuity of $F$ at $(A_1,A_2)$: if $x\in S(A_1,\epsilon_1),y\in S(A_2,\epsilon_1)$ then $|F(x,y)-F(A_1,A_2)|<\epsilon$. Then choose $\delta_i\ (i=1,2)$ such that $K\in[I]^{<\infty},\ d_H(K,I)<\delta_i$ implies that $|s_i(K)-A_i|<\epsilon_1$. Let $\delta=\min\{\delta_1,\delta_2\}$. If $K\in[I]^{<\infty},\ d_H(K,I)<\delta$ then $|s(K)-F(A_1,A_2)|<\epsilon$.
\end{proof}

\begin{cor}Let $s_1:[I]^{<\infty}\to\mathbb{R},\ s_2:[I]^{<\infty}\to\mathbb{R}$ and both ext$(s_1,I)$ and ext$(s_2,I)$ exist and be finite. Then ext$(\alpha_1s_1+\alpha_2s_2,I)\ \ (\alpha_1,\alpha_2\in\mathbb{R})$ exists and equals to $\alpha_1$ext$(s_1,I)+\alpha_2$ext$(s_2,I)$.
\end{cor}
\begin{proof}Apply \ref{plinearity} for $F(x_1,x_2)=\alpha_1x_1+\alpha_2x_2$.
\end{proof}

\begin{prp}If ext$(s,I)$ exists and finite then there is $\delta>0$ such that $\{s(K): d_H(K,I)<\delta\}$ is bounded.
\end{prp}
\begin{proof}Simply choose $\delta$ for $\epsilon=1$.
\end{proof}

\begin{prp}Let $s_n:[I]^{<\infty}\to\mathbb{R}\ (n\in\mathbb{N}),\ s:[I]^{<\infty}\to\mathbb{R}$. Suppose $s_n\to s$ on $[I]^{<\infty}$ uniformly in the following weaker sense: $\forall\epsilon>0\ \exists\delta>0\ \exists N\in\mathbb{N}$ such that $n>N, K\in[I]^{<\infty},d_H(K,I)<\delta$ implies that $|s_n(K)-s(K)|<\epsilon$. Let ext$(s_n,I)$ exists and be finite ($n\in\mathbb{N}$). Then $\lim\limits_{n\to\infty}$ext$(s_n,I)$ and ext$(s,I)$ exist and ext$(s,I)=\lim\limits_{n\to\infty}$ext$(s_n,I)$.
\end{prp}
\begin{proof}Let $A_n=$ext$(s_n,I)$. Take an accumulation point of $(A_n)$, say $A$. It is enough to prove that ext$(s,I)=A$.

First let $A\in\mathbb{R}$. Assume that $A_n\to A$. Let $\epsilon>0$. Choose $N_1\in\mathbb{N}$ such that $n>N_1$ implies that $|A_n-A|<\frac{\epsilon}{3}$. Choose $\delta_2>0,\ N_2\in\mathbb{N}$ such that $n>N_2, K\in[I]^{<\infty},d_H(K,I)<\delta_2$ implies that $|s_n(K)-s(K)|<\frac{\epsilon}{3}$. Let $N>\max\{N_1,N_2\}$. Choose $\delta_1>0$ such that $K\in[I]^{<\infty},d_H(K,I)<\delta_1$ implies that $|s_N(K)-$ext$(s_N,I)|<\frac{\epsilon}{3}$. Let $\delta=\min\{\delta_1,\delta_2\}$. Now if $K\in[I]^{<\infty},d_H(K,I)<\delta$ then we get that
\[|s(K)-A|\leq|s(K)-s_N(K)|+|s_N(K)-A_N|+|A_N-A|<\frac{\epsilon}{3}+\frac{\epsilon}{3}+\frac{\epsilon}{3}=\epsilon.\]

If $A=\pm\infty$ then a similar argument works..
\end{proof}

%--------------------------------
\subsection{Subspaces}

Now we examine relations of extensions to subspaces. 

\begin{thm}\label{tspacesextf}Let $I_1,I_2$ be subspaces of $I$ such that $d(I_1,I_2)>0$. Let $s_1=s||_{I_1},\ s_2=s||_{I_2}$. Suppose there is a continuous function $F:\mathbb{R}\times\mathbb{R}\to\mathbb{R}$ such that $s(K)=F(s_1(K\cap I_1),s_2(K\cap I_2))$ ($K\subset I$ finite). If ext$(s_1,I_1)$ and ext$(s_2,I_2)$ exist then so does ext$(s,I)$ and ext$(s,I)=F($ext$(s,I_1),$ext$(s,I_2))$.  
\end{thm}
\begin{proof}Let ext$(s_1,I_1)=A_1$ and ext$(s_2,I_2)=A_2$. Let $\epsilon>0$ and choose $\epsilon_1$ according to the continuity of $F$ at $(A_1,A_2)$: if $x\in S(A_1,\epsilon_1),y\in S(A_2,\epsilon_1)$ then $|F(x,y)-F(A_1,A_2)|<\epsilon$. Then choose $\delta_1,\delta_2$ such that $K_1\in[I_1]^{<\infty},\ d_H(K_1,I_1)<\delta_1$ implies that $|s_1(K_1)-A_1|<\epsilon_1$ and $K_2\in[I_2]^{<\infty},\ d_H(K_2,I_2)<\delta_2$ implies that $|s_2(K_2)-A_2|<\epsilon_1$. Let $\delta=\min\{d(I_1,I_2),\delta_1,\delta_2\}$. If $K\in[I]^{<\infty},\ d_H(K,I)<\delta$ then $d_H(K\cap I_1,I_1)<\delta$ and $d_H(K\cap I_2,I_2)<\delta$ hence 
$|s_1(K\cap I_1)-A_1|<\epsilon_1$ and $|s_2(K\cap I_2)-A_2|<\epsilon_1$. But $F(s_1(K\cap I_1),s_2(K\cap I_2))=s(K)$ therefore $|s(K)-F(A_1,A_2)|<\epsilon$.
\end{proof}

\begin{thm}\label{tdensesubsp}Let $I_1$ be a dense subspace of $I$ and $s_1=s||_{I_1}$. If ext$(s,I)$ exists then so does ext$(s_1,I_1)$ and ext$(s,I)$=ext$(s,I_1)$.  
\end{thm}
\begin{proof}If $K\subset I_1$ is finite and $\delta$-sdense in $I_1$ then $K$ is $\delta$-sdense in $I$.
\end{proof}

\begin{ex}Let $I_1,I_2$ be distinct dense subspaces of the interval $(0,1)$ and let ext$(s_1,I_1)$ and ext$(s_2,I_2)$ exist and ext$(s_1,I_1)=1$ and ext$(s_2,I_2)=2$. Let $s:[I]^{<\infty}\to\mathbb{R}$ such that $s||_{I_1}=s_1,\ s||_{I_2}=s_2$. Then ext$(s,I)$ does not exist by \ref{tdensesubsp}.
\end{ex}

\begin{ex}Let $I_1,I_2$ be distinct dense subspaces of the interval $I=(0,1)$ such that $I=I_1\cup I_2$. Let $s_1$ and $s_2$ be defined in the way that for sets containing an only point they take 1, for sets containing more than one point they take 0 and $s_1(\emptyset)=s_2(\emptyset)=0$. Obviously ext$(s_1,I_1)$ and ext$(s_2,I_2)$ exist and ext$(s_1,I_1)=$ext$(s_2,I_2)=0$. Let $s:[I]^{<\infty}\to\mathbb{R}$ such that 
\[s(K)=\frac{s_1(K\cap I_1)+s_2(K\cap I_2)}{2}\ \ (K\in[I]^{<\infty}).\]
Then ext$(s,I)$ does not exist because for each $\delta>0$ one can take a $\delta$-sdense set $D$ such that $D$ contains only one point from either $I_1$ or $I_2$ and then $s(D)=\frac{1}{2}$ and one can take a $\delta$-sdense set $E$ such that $E$ contains more than one point from both $I_1$ and $I_2$ and then $s(E)=0$.
\end{ex}

\begin{ex}Let $I=[-1,2]$ equipped with the usual metric. If $K\in [I]^{<\infty}$ then let 
\[s(K)=\begin{cases}\sum\limits_{k\in K}k&\text{if }K\subset[-1,1]\\
0&\text{otherwise.}
\end{cases}
\]
Clearly ext$(s,I)$ exists and equals to 0 while ext$(s||_{[-1,1]},[-1,1])$ does not exists.
\end{ex}

\begin{ex}Let $I=[-1,1]$ equipped with the usual metric. If $K\in [I]^{<\infty}$ then let 
$s(K)=\sum\limits_{k\in K}k$.
Clearly ext$(s,I)$ does not exists however it is easy to see that for every $x\in\overline{\mathbb{R}}$ there is $H\subset I$ such that ext$(s||_H,H)=x$.
\end{ex}

%--------------------------------
\subsection{Increasing/decreasing related properties}

We now look for conditions for the extension being equal to the supremum of values of $s$.

\begin{df}A function $s:[I]^{<\infty}\to\mathbb{R}$ is called \textbf{increasing} (decreasing) if $K,L\in[I]^{<\infty},\ K\subset L$ implies that $s(K)\leq s(L)$\ \ ($s(K)\geq s(L)$).
\end{df}

\begin{thm}\label{tsincextsup}Let $s$ be increasing (decreasing). If ext$(s,I)$ exists then it equals to $\sup\{s(K):K\in [I]^{<\infty}\}$\ \ ($\inf\{s(K):K\in [I]^{<\infty}\}$).
\end{thm}
\begin{proof}Let $A=$ext$(s,I),\ r=\sup\{s(K):K\in [I]^{<\infty}\},\ r\in\mathbb{R}$. Obviously $r<A$ cannot happen. Suppose that $A<r$. Take a set $L\in [I]^{<\infty}$ such that $s(L)>\frac{A+r}{2}$. Then for $\epsilon=\frac{r-A}{2}$ find $\delta$ according to the definition of ext$(s,I)$. If we take a set $K\in [I]^{<\infty}$ that is $\delta$-sdense then so is $K\cup L$. But by $s$ being increasing $s(K\cup L)\geq s(L)>\frac{A+r}{2}$ which is a contradiction.

If $r=\infty$ then a similar argument works.

The decreasing case can be handled similarly.
\end{proof}

\begin{ex}\label{eincnotdinc}Let $I_0,I_1$ be distinct dense subspaces of the interval $I=(0,1)$ such that $I=I_0\cup I_1$. For $K\in[I]^{<\infty}$ let 
\[s(K)=\begin{cases}
0&\text{if }K\subset I_0\\
1&\text{otherwise.}
\end{cases}\]
Obviously $s$ is increasing but ext$(s,I)$ does not exist because if $K_i\in[I_i]^{<\infty}$ is $\delta$-sdense set on $I_i\ (i=0,1)$ then it is $\delta$-sdense set on $I$ as well and $s(K_0)=0$ while $s(K_1)=1$.
\end{ex}

\begin{df}\label{ddinc}A function $s:[I]^{<\infty}\to\mathbb{R}$ is called \textbf{d-increasing} (d-decreasing) if $\forall\epsilon>0$ and for all $K\in[I]^{<\infty}$ there is $\delta>0$ such that $L\in[I]^{<\infty}$ being $\delta$-sdense implies that $s(K)-\epsilon<s(L)$\ \ ($s(K)+\epsilon>s(L)$).
\end{df}

\begin{thm}\label{pdincsup}If $s:[I]^{<\infty}\to\mathbb{R}$ be d-increasing (d-decreasing) then ext$(s,I)$ exists and equals to $\sup\{s(K):K\in [I]^{<\infty}\}$\ \ ($\inf\{s(K):K\in [I]^{<\infty}\}$).
\end{thm}
\begin{proof}Let $S=\sup\{s(K):K\in [I]^{<\infty}\}\in\mathbb{R}$ and $\epsilon>0$. Choose $K\in [I]^{<\infty}$ such that $S-\frac{\epsilon}{2}<s(K)$ and choose $\delta>0$ such that $L\in[I]^{<\infty}$ being $\delta$-sdense implies that $s(K)-\frac{\epsilon}{2}<s(L)$. Then $S-\epsilon<s(L)$.

Let $S=+\infty$ and $l\in\mathbb{R}$. Choose $K\in [I]^{<\infty}$ such that $l+1<s(K)$ and choose $\delta>0$ such that $L\in[I]^{<\infty}$ being $\delta$-sdense implies that $s(K)-1<s(L)$. Then $l<s(L)$.

The case when $S=-\infty$ can be handled similarly.
\end{proof}

\begin{prp}If ext$(s,I)=+\infty\ (-\infty)$ then $s$ is d-increasing (d-decreasing).\qed
\end{prp}

\begin{ex}\label{einnJordm}Let $H\subset\mathbb{R}$ be bounded. Let $I=H$. If $K\in[I]^{<\infty},K=\{k_1,\dots,k_m\},k_1<\dots<k_m$ then set 
\[s(K)=\sum\limits_{\substack{1\leq i\leq m-1\\(k_i,k_{i+1})\subset H}}k_{i+1}-k_i.\]
We show that $s$ is d-increasing. Let $\epsilon>0$. Let $L\in[I]^{<\infty},\ d_H(L,I)<\frac{\epsilon}{2m}$. If $1\leq i\leq m-1,(k_i,k_{i+1})\subset H$ then it is obvious that those intervals which endpoints $\in L\cap(k_i,k_{i+1})$ cover at least $k_{i+1}-k_i-\frac{\epsilon}{m}$ length from $(k_i,k_{i+1})$. Therefore it gives that $s(K)-\epsilon<s(L)$.

By \ref{pdincsup} we get that ext$(s,I)$ is the inner Jordan measure of $H$.

\smallskip

Note that $s$ is increasing as well.
\end{ex}

\begin{ex}None of the properties increasing and d-increasing implies the other.

We have already noted that example in \ref{eincnotdinc} is increasing. It is not d-increasing since let $K$ a singleton from $I_1$. Then $s(K)=1$ however for all $\delta>0$ we can take $L\in[I_0]^{<\infty}$ such that $d_H(L,I)<\delta$. But $s(L)=0$.

\smallskip

Let a series be given such that $\sum\limits_{i=1}^{\infty} a_i=+\infty,\ \sum\limits_{i=1}^{\infty} |a_i|^-<\infty$ and $\{i\in\mathbb{N}:a_i<0\}$ is infinite. Let $s,I$ be defined as in subsection \ref{ssabsconvser}. It can be readily seen that $s$ is d-increasing however fails to be increasing.

\smallskip

One can also create an example where ext$(s,I)$ is finite, $s$ is d-increasing but not increasing. E.g. let $I=\left\{\frac{1}{2^n}:n\in\mathbb{N}\right\}$ equipped with the usual metric, for $K\in[I]^{<\infty}$ let
\[s(K)=
\begin{cases}
\sum\limits_{k\in K}k&\text{if }|K|\text{ is even}\\
\sum\limits_{\substack{k\in K\\ k\ne\min K}}k-\frac{1}{|K|}&\text{if }|K|\text{ is odd.}
\end{cases}\]
It is easy to see that all required conditions hold. Moreover $\forall\delta>0\ \exists K\in[I]^{<\infty}$ such that $d_H(K,I)<\delta$ and $\exists L\in[I]^{<\infty}$ such that $K\subset L$ but $s(K)\not\leq s(L)$.
\end{ex}

%--------------------------------
\subsection{On continuity}

\begin{df}A function $s:[I]^{<\infty}\to\mathbb{R}$ is called \textbf{d-continuous} if $\forall n\in\mathbb{N}\ s|_{[I]^n}$ is continuous regarding the Hausdorff pseudo-metric on $[I]^n$.
\end{df}

When $I$ is ordered then $[I]^n$ can be identified with a subset of $I^n$: If $K\in[I]^n,\ K=\{k_1,\dots,k_n\},\ k_1<\dots<k_n$ then $f(K)=(k_1,\dots,k_n)\in I^n$ is an injective mapping.

\begin{prp}Let $I$ is ordered and $s:[I]^{<\infty}\to\mathbb{R}$. Then $s$ is d-continuous iff $\forall n\in\mathbb{N}\ s|_{[I]^n}$ is continuous considering $s|_{[I]^n}$ as a $Q\to\mathbb{R}$ function with the identification described above where $Q=f([I]^n)$.
\end{prp}
\begin{proof}If $s$ is d-continuous, $n\in\mathbb{N},\ K\in[I]^n,K=\{k_1,\dots,k_n\},k_1<\dots<k_n,\ \epsilon>0$, choose $\delta$ according to the definition. If $L=\{l_1,\dots,l_n\},l_1<\dots<l_n$ and $d(f(K),f(L))<\delta$ then $d(k_i,l_i)<\delta\ (1\leq i\leq n)$ which gives that $d_H(K,L)<\delta$ hence $|s(K)-s(L)|<\epsilon$.

If $s|_{[I]^n}:Q\to\mathbb{R}$ is continuous, $K\in[I]^n,K=\{k_1,\dots,k_n\},k_1<\dots<k_n,\ \epsilon>0$ then choose $\delta$ according to the definition. Let $\delta_1<\delta$ and $\delta_1<\frac{1}{2}\min\{d(x,y):x,y\in K\}$; set $\delta_2=\frac{\delta_1}{\sqrt{n}}$. Now if $L\in[I]^n,L=\{l_1,\dots,l_n\},l_1<\dots<l_n$,  $d_H(K,L)<\delta_2$ then $d(k_i,l_i)<\delta_2\ (1\leq i\leq n)$ hence $d(f(K),f(L))<\delta$ which gives that $|s(K)-s(L)|<\epsilon$.
\end{proof}

\begin{prp}If $s:[I]^{<\infty}\to\mathbb{R}$ is continuous regarding the Hausdorff pseudo-metric on $[I]^{<\infty}$ then $s$ is d-continuous.\qed
\end{prp}

The continuity of $s$ is a very strong condition as straightforward examples can demonstrate.

\begin{ex}In general neither finite summation nor finite average is continuous regarding the Hausdorff pseudo-metric on finite sets.
\end{ex}
\begin{proof}For finite summation let $I=\{0\}\cup\{\frac{1}{n}:n\in\mathbb{N}\}$ and set $s(K)=\sum\limits_{k\in K}k\ (K\in[I]^{<\infty})$. Let $K=\{0\}$. Clearly for each $\delta>0$ there is $L\in[I]^{<\infty}$ such that $d_H(K,L)<\delta$ and $s(L)>1$ because $\sum\limits_{n=1}^{\infty}\frac{1}{n}=\infty$.

For finite average let $I=\{0,2\}\cup\{\frac{1}{n}:n\in\mathbb{N}\}\cup\{2-\frac{1}{n}:n\in\mathbb{N}\}$ and set $s(K)=\A(K)$. Let $K=\{0,2\}$. Then $s(K)=1$. Clearly for each $\delta>0$ there are $L_1,L_2\in[I]^{<\infty}$ such that $d_H(K,L_i)<\delta\ (i=1,2)$ and $s(L_1)<0.5,\ s(L_2)>1.5$.
\end{proof}

Here we present a few examples where $s$ is continuous.

\begin{ex}If $K\in[I]^{<\infty}$ then set $s(K)=\frac{\min K+\max K}{2}$. Clearly $s$ is continuous on $[I]^{<\infty}$ and ext$(s,I)=\frac{\inf I+\sup I}{2}$.
\end{ex}

\begin{ex}If $K\in[I]^{<\infty}$ then set $s(K)=\mathrm{diam}\ K=\max\{d(x,y):x,y\in K\}$. Clearly $s$ is continuous on $[I]^{<\infty}$ and ext$(s,I)=\mathrm{diam}\ I$.
\end{ex}

\begin{prp}\label{parclength}Let $I=[0,1],\ \gamma:I\to\mathbb{R}^n$ continuous and rectifiable. If for $K\in[I]^{<\infty}$ $s(K)$ is the length of the inscribed polygon determined by $\gamma$ and the points of $K$ then $s$ is continuous on $[I]^{<\infty}$.
\end{prp}
\begin{proof}Let $l$ be the length function on $I$ i.e. $l:I\to\mathbb{R}, \ x\in I,\ l(x)$ is the length of the curve $\gamma|_{[0,x]}$. It is known that $l$ is continuous assuming that $\gamma$ is rectifiable.

Let $\epsilon>0$ and $K\in[I]^{<\infty}, m=|K|, K=\{k_1,\dots,k_m\}, k_1<\dots<k_m$. Then there is $\delta_1>0$ such that $\forall k\in K$ the length of $\gamma|_{[k-\delta_1,k+\delta_1]}<\frac{\epsilon}{2m}$. Using the continuity of $\gamma$ there is $\delta_2>0$ such that $\forall i (1\leq i\leq m-1)\ \forall k_i'\in I\ |k_i-k_i'|<\delta_2$ implies that the difference between the length of the line joining $\gamma(k_i)$ and $\gamma(k_{i+1})$ and the length of the line joining $\gamma(k_i')$ and $\gamma(k'_{i+1})$ is less than $\frac{\epsilon}{2m}$.

Let $\delta<\delta_1,\delta<\delta_2,\delta<\frac{1}{2}\min\{|x-y|:x,y\in K,x\ne y\}$. Now if $d_H(K,K')<\delta$ then it is easy to see that $|s(K)-s(K')|<\epsilon$.
\end{proof}

\begin{prp}If $s:[I]^{<\infty}\to\mathbb{R}$ is increasing (decreasing) and d-continuous then $s$ is d-increasing (d-decreasing).
\end{prp}
\begin{proof}Let $K\in[I]^{n}$ and $\epsilon>0$. Choose $\delta$ according to the d-continuity for $n$. Let $\delta_1<\delta$ and $\delta_1<\min\{d(x,y):x,y\in K\}$. If $L\in[I]^{<\infty}$ is $\frac{\delta_1}{2}$-sdense then for each $k\in K$ choose distinct $l_k\in L$ such that $d(k,l_k)<\delta_1$. It can be done. Let $L_1=\{l_k:k\in K\}$. Then $L_1\in[I]^{n},\ d_H(K,L_1)<\delta$ hence $|s(K)-s(L_1)|<\epsilon$. $L_1\subset L$ implies that $s(L_1)\leq s(L)$ which altogether gives that $s(K)-\epsilon<s(L)$.

The decreasing part is similar.
\end{proof}

\begin{ex}The $s$ in example \ref{einnJordm} is increasing, d-increasing however not d-continuous. To show that let $I=[0,1]\cup([1,2]\cap\mathbb{Q}),\ K=\{0,1\}$. Then $s(K)=1$ but there is $L\in[I]^2$ arbitrarily close to $K$ such that $s(L)=0$.
\end{ex}

\begin{cor}Both finite summation on positive numbers and length of finite polygons are d-increasing.
\end{cor}
\begin{proof}It is easy to see that they are increasing and d-continuous.
\end{proof}

\begin{cor}\label{carclendefeq}The two definitions of arc length are equivalent (see section \ref{sbnn}).
\end{cor}
\begin{proof}The second definition gives exactly ext$(s,I)$ using the notations from \ref{parclength}.
\end{proof}

\begin{ex}\label{edarboux}If a bounded $f:[a,b]\to\mathbb{R}$ is given then let $I=[a,b]$ equipped with the usual metric and $s:[I]^{<\infty}\to\mathbb{R}$ defined in the following way. Let $H\subset I$ be finite. Let $\tilde H=H\cup\{a,b\}$. Let $\tilde H=\{a_0,\dots,a_n\}$ where $a=a_0<\dots <a_n=b$. Now set
\[s(H)=\sum\limits_{i=0}^{n-1}\sup f\big([a_i,a_{i+1}]\big)\cdot (a_{i+1}-a_i)\]
which is the upper (Darboux) sum corresponding to the partition of $\tilde H$. Clearly $s$ is decreasing. It is easy to see that ext$(s,I)$ exists and equals to the upper Darboux integral of $f$ (actually the two definitions are almost identical). Hence by \ref{tsincextsup} we get a new proof for the well known fact that the upper Darboux integral equals to $\inf\{s(K):K\in [I]^{<\infty}\}$.

We show that $s$ is not d-decreasing in general. In order to see that let $I=[0,1],\ H=\{0,1\},$
\[f(x)=
\begin{cases}
0&\text{if }x\ne 0.5\\
1&\text{if }x=0.5.
\end{cases}
\]
Then $s(H)=1$. If $\delta>0$ then let $\delta'=\min\{\delta,0.1\}$ and $L=\{i\delta':i\in\mathbb{Z},\ a\leq i\delta'\leq b\}$. Then clearly $L\in[I]^{<\infty}$, $\delta$-sdense and  $s(L)<0.1$.  

The same example shows that $s$ is not d-continuous in general, just take $H=\{0,0.5,1\}$. If we move 0.5 slightly then $s$ changes by 0.25 at least.
\end{ex}

\begin{ex}Very similar arguments and the same example show that $s$ defined in subsection \ref{ssriemann} (for Riemann integral) is neither d-decreasing nor d-continuous in general.

It is also not true that for continuous $f$, $s$ would be continuous. Let
\[f(x)=
\begin{cases}
0&\text{if }0\leq x\leq 0.25\\
4x-1&\text{if }0.25\leq x\leq 0.5\\
3-4x&\text{if }0.5\leq x\leq 0.75\\
0&\text{if }0.75\leq x\leq 1.
\end{cases}\]
Let $H=\{0,0.5,1\},\ \delta<0.1$. It is easy to choose $L\in[I]^{<\infty}$ such that $d_H(L,H)<\delta,\ L=\{l_0,\dots,l_n\},\ l_0<\dots<l_n$ and there is $i$ such that the contribution of the three points $l_i,l_{i+1},l_{i+2}$ to $s(L)$ is $f(l_{i+1})(l_{i+2}-l_{i})$ and $l_{i+1}<0.1$ while $l_{i+2}>0.4$. Then $f(l_{i+1})(l_{i+2}-l_{i})=0$ which gives that $s(L)<0.6$ while $s(H)=1$.

On the original domain $s$ is not d-continuous either. Take the same $f$ and $H$. Let $K=\{\epsilon,0.5,1\}$ for $\epsilon<0.25$. Then $\tilde K=\{0,\epsilon,0.5,1\}$, $|K|-1$ is odd hence $s(K)=f(\epsilon)(0.5-0)+f(1)(1-0.5)=0\cdot 0.5+0\cdot 0.5=0$ (see \ref{ssriemann}).
\end{ex}

However on a slightly modified domain d-continuity follows:

\begin{thm}If $f:[a,b]\to\mathbb{R}$ is continuous, $s:\big[(a,b)\big]^{<\infty}\to\mathbb{R}$ and $s$ is defined as in subsection \ref{ssriemann} then $s$ is d-continuous.
\end{thm}
\begin{proof}Suppose that $|f|<M$. Let $\epsilon>0,\ H\subset (a,b)$ be finite, $H=\{a_1,\dots,a_n\},\ a_1<\dots <a_n$. Let $\tilde H=H\cup\{a,b\}$. For $\frac{\epsilon}{2(b-a)}$ choose $\delta_1$ according to the uniform continuity of $f$ and let 
\[\delta<\frac{1}{2}\min\big\{|x-y|:x,y\in\tilde H,x\ne y\big\}\text{ and }\delta<\min\left\{\delta_1,\frac{\epsilon}{4nM}\right\}.\]
Take $K\in\big[(a,b)\big]^{<\infty}$ such that $K=\{b_1,\dots,b_n\},\ |a_i-b_i|<\delta\ (1\leq i\leq n)$. Let $a_0=b_0=a,a_{n+1}=b_{n+1}=b$. Then there are two cases depending on the parity of $n$. If $n=2k$ then
\[\big|s(H)-s(K)\big|=\left|\sum\limits_{i=0}^{k-1}f(a_{2i+1})(a_{2i+2}-a_{2i})-\sum\limits_{i=0}^{k-1}f(b_{2i+1})(b_{2i+2}-b_{2i})\right|\leq\]
\[\sum\limits_{i=0}^{k-1}\big|f(a_{2i+1})(a_{2i+2}-a_{2i})-f(b_{2i+1})(a_{2i+2}-a_{2i})\big|+\big|f(b_{2i+1})(a_{2i+2}-a_{2i})-f(b_{2i+1})(b_{2i+2}-b_{2i})\big|\leq\]
\[\sum\limits_{i=0}^{k-1}\big|f(a_{2i+1})-f(b_{2i+1})\big|\big|a_{2i+2}-a_{2i}\big|+\sum\limits_{i=0}^{k-1}\big|f(b_{2i+1})\big|\big|(a_{2i+2}-a_{2i})-(b_{2i+2}-b_{2i})\big|\leq\]
\[\frac{\epsilon}{2(b-a)}(b-a)+M\sum\limits_{i=0}^{k-1}|a_{2i+2}-b_{2i+2}|+|a_{2i}-b_{2i}|\leq\frac{\epsilon}{2}+2kM\frac{\epsilon}{4nM}=\epsilon.\]

If $n=2k+1$ then a similar argument works.
\end{proof}

%--------------------------------
\subsection{Integral on measure spaces}

Now we show that the integral on measure spaces can be fit into this framework too.

\begin{prp}\label{pintinmeassp}Let $\langle X,\mu\rangle$ be a measure space such that $\mu(X)<\infty$. Let $f:X\to\mathbb{R}$ be a measurable function such that $0\leq f<M\in\mathbb{R}$. Let $I=(0,M)$ and $s:[I]^{<\infty}\to\mathbb{R}$ defined as follows. If $H\in[I]^{<\infty},\ H=\{a_1,\dots,a_n\},\ a_1<\dots<a_n$ then let $a_0=0,a_{n+1}=M$ and set
\[s(H)=\sum\limits_{i=0}^{n}a_i\mu\Big(f^{-1}\big([a_i,a_{i+1})\big)\Big).\]
Then $s$ is d-increasing.
\end{prp}
\begin{proof}
Let $H\in[I]^{<\infty},\ \epsilon>0$. Let $\delta<\frac{\epsilon}{\mu(X)}$. Let $K\in[I]^{<\infty},\ d_H(K,I)<\delta$. Let $H\cup K=\{c_1,\dots,c_l\},\ c_1<\dots<c_l,\ c_0=0,c_{l+1}=M$. Let 
\[h:H\cup K\to H\cup\{0\},\ h(c_i)=\max\{x\in H\cup\{0\}:x\leq c_i\},\]
\[k:H\cup K\to K\cup\{0\},\ k(c_i)=\max\{x\in K\cup\{0\}:x\leq c_i\}.\]
Now we get that
\[s(H)-s(K)=\sum\limits_{i=0}^{l}h(c_i)\mu\Big(f^{-1}\big([c_i,c_{i+1})\big)\Big)-\sum\limits_{i=0}^{l}k(c_i)\mu\Big(f^{-1}\big([c_i,c_{i+1})\big)\Big)=\]
\[\sum\limits_{i=0}^{l}\big(h(c_i)-k(c_i)\big)\mu\Big(f^{-1}\big([c_i,c_{i+1})\big)\Big)\leq\]
\[\sum\limits_{h(c_i)>k(c_i)}\big(h(c_i)-k(c_i)\big)\mu\Big(f^{-1}\big([c_i,c_{i+1})\big)\Big)\leq\delta\mu(X)<\epsilon\]
because if $h(c_i)>k(c_i)$ then $c_i\in K$ cannot happen, and if $c_i\in H$ then $h(c_i)-k(c_i)<\delta$.
\end{proof}

\begin{thm}\label{tintmeassp}Let $\langle X,\mu\rangle$ be a measure space such that $\mu(X)<\infty$. Let $f:X\to\mathbb{R}$ be a measurable function such that $0\leq f<M\in\mathbb{R}$. Let $I=(0,M)$ and $s:[I]^{<\infty}\to\mathbb{R}$ defined as in \ref{pintinmeassp}. Then

(a) ext$(s,I)$ exists and equals to $\int\limits_X f d\mu$.

(b) For a given function $f$ add $f$ to the notation of $s$: $s_f$. If $f$ can take both positive and negative values then ext$(s_{|f|^+},I)-$ext$(s_{|f|^-},I)=\int\limits_X f d\mu$.
\end{thm}
\begin{proof}
(a) By \ref{pdincsup} ext$(s,I)$ exists and equals to $\sup\{s(K):K\in [I]^{<\infty}\}$ which is known to equal to $\int\limits_X f d\mu$ (see \cite{bruckner2}).

(b) It is a straightforward consequence of (a).
\end{proof}

We show that the integral can also be represented as one ext.

The following lemma is straightforward hence the proof is omitted.

\begin{lem}If $s_1:[I]^{<\infty}\to\mathbb{R}$ is d-increasing and $s_2:[I]^{<\infty}\to\mathbb{R}$ is d-decreasing then $s_1-s_2$ is d-increasing.\qed
\end{lem}

\begin{thm}Let $\langle X,\mu\rangle$ be a measure space such that $\mu(X)<\infty$. Let $f:X\to\mathbb{R}$ be a measurable function such that $|f|<M\in\mathbb{R}$. Let $I=(-M,M)-\{0\}$ and $s:[I]^{<\infty}\to\mathbb{R}$ defined as follows. If $H\in[I]^{<\infty},\ H\cup\{0\}=\{a_1,\dots,a_n\},\ a_1<\dots<a_n$ then let $a_0=-M,a_{n+1}=M$ and set
\[s(H)=\sum\limits_{i=0}^{n}a_i\mu\Big(f^{-1}\big([a_i,a_{i+1})\big)\Big).\]
Then ext$(s,I)$ exists and equals to $\int\limits_X f d\mu$.
\end{thm}
\begin{proof} If $H\in[I]^{<\infty},\ H\cup\{0\}=\{a_1,\dots,a_n\},\ a_1<\dots<a_n$ then let $a_0=-M,a_{n+1}=M$ and set 
\[s_1(H)=\sum\limits_{\substack{{0\leq i\leq n}\\{0\leq a_i}}}a_i\mu\Big(f^{-1}\big([a_i,a_{i+1})\big)\Big),\ 
s_2(H)=\sum\limits_{\substack{{0\leq i\leq n}\\{a_i<0}}}a_i\mu\Big(f^{-1}\big([a_i,a_{i+1})\big)\Big).\]
Clearly $s=s_1+s_2$. As in \ref{pintinmeassp} one can show that $s_1$ is d-increasing while $-s_2$ is d-decreasing. Hence ext$(s_1,I)$ exists and equals to $\sup\{s_1(K):K\in [I]^{<\infty}\}$ which equals to $\int\limits_X |f|^+ d\mu$. Similarly ext$(-s_2,I)$ exists and equals to $\inf\{-s_2(K):K\in [I]^{<\infty}\}$ which equals to $\int\limits_X |f|^- d\mu$. Also $s$ is d-increasing so ext$(s,I)$ exists and equals to $\sup\{s(K):K\in [I]^{<\infty}\}=\sup\{s_1(K):K\in [I]^{<\infty}\}-\inf\{-s_2(K):K\in [I]^{<\infty}\}=\int\limits_X |f|^+ d\mu - \int\limits_X |f|^- d\mu=\int\limits_X f d\mu$.
\end{proof}

Using the simplest form we show some more properties.

\begin{prp}Using the notations of \ref{pintinmeassp} it can be readily seen that $s$ is increasing as well.\qed
\end{prp}

\begin{ex}\label{eaintnotdcont}Using the notations of \ref{pintinmeassp} we show that $s$ is not continuous moreover not even d-continuous in general. To show that let $X=[0,1]$ equipped with the Lebesgue measure, $f:[0,1]\to\mathbb{R},$
\[f(x)=\begin{cases}
1&\text{if }x\ne 1\\
2&\text{if }x=1\\
\end{cases}\]
and let $H=\{1\}$. Then $s(H)=1$ however if $K=\{k\}$ with $k>1$ then $s(K)=0$.
\end{ex}

However we show that $s$ is left-continuous in the following sense.

\begin{df}A function $s:[I]^{<\infty}\to\mathbb{R}$ is called \textbf{left-continuous} if $\forall\epsilon>0$ for all $H\in[I]^{<\infty}\ \exists\delta>0$ such that $K\in[I]^{<\infty},\ d_H(H,K)<\delta,\ \forall b\in K\ \exists a\in H$ such that $b\leq a,\ |a-b|<\delta$ implies that $|s(H)-s(K)|<\epsilon$.
\end{df}

\begin{df}A function $s:[I]^{<\infty}\to\mathbb{R}$ is called \textbf{left-d-continuous} if $\forall\epsilon>0$ for all $H\in[I]^{<\infty}\ \exists\delta>0$ such that $K\in[I]^{<\infty},\ |K|=|H|, H=\{a_1,\dots,a_n\},a_1<\dots<a_n, K=\{b_1,\dots,b_n\},b_1<\dots<b_n,\ \forall i\in\{1,\dots,n\},\ b_i\leq a_i,\ |a_i-b_i|<\delta$ implies that $|s(H)-s(K)|<\epsilon$.
\end{df}

\begin{rem}If $s:[I]^{<\infty}\to\mathbb{R}$ is left-continuous then it is clearly left-d-continuous as well.
\end{rem}

\begin{thm}Let $\langle X,\mu\rangle$ be a measure space such that $\mu(X)<\infty$. Let $f:X\to\mathbb{R}$ be a measurable function such that $0\leq f<M\in\mathbb{R}$. Let $I=(0,M)$ and $s:[I]^{<\infty}\to\mathbb{R}$ defined as in \ref{pintinmeassp}. Then $s$ is left-continuous.
\end{thm}
\begin{proof}Let $\epsilon>0,\ H\in[I]^{<\infty},\ H=\{a_1,\dots,a_n\},a_1<\dots<a_n$, let $a_0=0,a_{n+1}=M$. Let $\delta_1=\frac{\epsilon}{2\mu(X)},\ \delta_2=\frac{1}{2}\min\{|x-y|:x,y\in H\cup\{0\},x\ne y\}$. For every $i\in\{1,\dots,n\}$
\[\bigcap\limits_{j=1}^{\infty}\mu\left(f^{-1}\left(\left[a_i-\frac{1}{j},a_i\right)\right)\right)=\emptyset\]
because $\bigcap\limits_{j=1}^{\infty}f^{-1}\left(\left[a_i-\frac{1}{j},a_i\right)\right)=\emptyset$ and $\mu\left(f^{-1}\left(\left[0,a_i\right)\right)\right)<+\infty$. Hence for every $i\in\{1,\dots,n\}$ there is $\gamma_i>0$ such that $b\leq a_i, |a_i-b|<\gamma_i$ implies that $\mu\left(f^{-1}\left(\left[b,a_i\right)\right)\right)<\frac{\epsilon}{8nM}$.

Now set $\delta=\min\{\delta_1,\delta_2,\gamma_i:1\leq i\leq n\}$.

\smallskip

Let $K\in[I]^{<\infty},\ d_H(H,K)<\delta,\ \forall b\in K\ \exists a\in H$ such that $b\leq a,\ |a-b|<\delta$. Let $K=\{b_1,\dots,b_m\},b_1<\dots<b_m$. Set $b_0=0,b_{m+1}=M$. Set $j_i=\max\{j:0\leq j\leq m+1,b_j\leq a_i\}\ \ (0\leq i\leq n+1)$. Obviously $b_{j_i}\leq a_i$ and $|a_i-b_{j_i}|<\delta$. Moreover if $j_i+1\leq k\leq j_{i+1}$ then $|a_{i+1}-b_{k}|<\delta$.

Then we get that
\[|s(H)-s(K)|=\left|\sum\limits_{i=0}^{n}a_i\mu\Big(f^{-1}\big([a_i,a_{i+1})\big)\Big)-\sum\limits_{j=0}^{m}b_i\mu\Big(f^{-1}\big([b_i,b_{i+1})\big)\Big)\right|=\]
\[\left|\sum\limits_{i=0}^{n}a_i\mu\Big(f^{-1}\big([a_i,a_{i+1})\big)\Big)-\sum\limits_{i=0}^{n}\sum\limits_{k=j_i}^{j_{i+1}-1}b_k\mu\Big(f^{-1}\big([b_k,b_{k+1})\big)\Big)\right|\leq\]
\[\sum\limits_{i=0}^{n}\left|a_i\mu\Big(f^{-1}\big([a_i,a_{i+1})\big)\Big)-b_{j_i}\mu\Big(f^{-1}\big([b_{j_i},b_{j_i+1})\big)\Big)-\sum\limits_{k=j_i+1}^{j_{i+1}-1}b_k\mu\Big(f^{-1}\big([b_k,b_{k+1})\big)\Big)\right|\leq\]
\[\sum\limits_{i=0}^{n}\left|a_i\mu\Big(f^{-1}\big([a_i,a_{i+1})\big)\Big)-b_{j_i}\mu\Big(f^{-1}\big([b_{j_i},b_{j_i+1})\big)\Big)\right|+\left|\sum\limits_{k=j_i+1}^{j_{i+1}-1}b_k\mu\Big(f^{-1}\big([b_k,b_{k+1})\big)\Big)\right|\leq\]
\[\sum\limits_{i=0}^{n}\left|a_i\mu\Big(f^{-1}\big([a_i,a_{i+1})\big)\Big)-b_{j_i}\mu\Big(f^{-1}\big([b_{j_i},b_{j_i+1})\big)\Big)\right|+\sum\limits_{i=0}^{n}\sum\limits_{k=j_i+1}^{j_{i+1}-1}\left|b_k\right|\mu\Big(f^{-1}\big([b_k,b_{k+1})\big)\Big).\]
Now we can estimate the first term as
\[\sum\limits_{i=0}^{n}\left|a_i\mu\Big(f^{-1}\big([a_i,a_{i+1})\big)\Big)-b_{j_i}\mu\Big(f^{-1}\big([b_{j_i},b_{j_i+1})\big)\Big)\right|\leq\]
\[\sum\limits_{i=0}^{n}\left|a_i\mu\Big(f^{-1}\big([a_i,a_{i+1})\big)\Big)-a_i\mu\Big(f^{-1}\big([b_{j_i},b_{j_i+1})\big)\Big)\right|+\]
\[\left|a_i\mu\Big(f^{-1}\big([b_{j_i},b_{j_i+1})\big)\Big)-b_{j_i}\mu\Big(f^{-1}\big([b_{j_i},b_{j_i+1})\big)\Big)\right|\leq\]
\[\sum\limits_{i=0}^{n}|a_i|\left(\mu\Big(f^{-1}\big([b_{j_i+1},a_{i+1})\big)\Big)+\mu\Big(f^{-1}\big([b_{j_i},a_i)\big)\Big)\right)+\sum\limits_{i=0}^{n}|a_i-b_{j_i}|\mu\Big(f^{-1}\big([b_{j_i},b_{j_i+1})\big)\Big)\leq\]
\[M\cdot n\cdot 2\cdot\frac{\epsilon}{8nM}+\delta\mu(X)=\frac{\epsilon}{4}+\frac{\epsilon}{2}=\frac{3\epsilon}{4}.\]
We can estimate the second term as
\[\sum\limits_{i=0}^{n}\sum\limits_{k=j_i+1}^{j_{i+1}-1}\left|b_k\right|\mu\Big(f^{-1}\big([b_k,b_{k+1})\big)\Big)\leq
\sum\limits_{i=0}^{n}M\mu\Big(f^{-1}\big([b_{j_i+1},b_{j_{i+1}})\big)\Big)\leq\]
\[\sum\limits_{i=0}^{n}M\mu\Big(f^{-1}\big([b_{j_i+1},a_{i+1})\big)\Big)\leq M\cdot n\cdot\frac{\epsilon}{8nM}=\frac{\epsilon}{8}.\]
Therefore we get that $|s(H)-s(K)|<\epsilon$.
\end{proof}

\begin{ex}We can examine the following weaker left-continuity type notion: $\forall\epsilon>0$ for all $H\in[I]^{<\infty}\ \exists\delta>0$ such that $K\in[I]^{<\infty},\ d_H(H,K)<\delta,\ \forall a\in H\ \exists b\in K$ such that $b\leq a,\ |a-b|<\delta$ would imply that $|s(H)-s(K)|<\epsilon$.

However it does not hold in general for $s$ defined in \ref{pintinmeassp}. Take the example in \ref{eaintnotdcont}. Let $\epsilon=0.5,\ K=\{0.9,1.1\}$. Then $s(K)=0+0.2\cdot 1+0=0.2$.
\end{ex}

%--------------------------------
\subsection{Approximating from outside}

We now examine if the set $I$ is approximating from outside too (from outside=from a bigger set) then when we get the same extension.

\begin{df}Let $\langle I,d\rangle$ be a totally bounded pseudo-metric space. Let $s:[I]^{<\infty}\to\mathbb{R}$. Let $J$ is a subspace of $I$. The \textbf{extension of $\boldsymbol s$ onto $\boldsymbol J$ within $\boldsymbol I$} equals to $A$ if $\forall\epsilon>0\ \exists\delta>0$ such that $K\in[I]^{<\infty}$ and $d_H(K,J)<\delta$ implies that $|s(K)-A|<\epsilon$. In this case we use the notation ext$(s,J,I)=A$.
\end{df}

\begin{prp}If $J$ is a subspace of $I$ and ext$(s,J,I)$ exists then so does ext$(s,J)$ and they are equal.
\end{prp}
\begin{proof}If $K\in[J]^{<\infty},\ d_H(K,J)<\delta$ holds in $J$ then $K\in[I]^{<\infty}$ and $d_H(K,J)<\delta$ holds in $I$.
\end{proof}

\begin{prp}\label{pdensesp2exteq}Let $J$ be a dense subspace of $I$. Then ext$(s,J,I)$ exists iff ext$(s,I)$ exists and they are equal.
\end{prp}
\begin{proof}If ext$(s,J,I)$ exists then note that $K\in[I]^{<\infty},\ d_H(K,I)<\delta$ implies that $d_H(K,J)<\delta$.

If ext$(s,I)$ exist then observe that $K\in[I]^{<\infty},\ d_H(K,J)<\delta$ implies that $d_H(K,I)<2\delta$.
\end{proof}

Now we present a simpler notion than the overall continuity of $s$ that will suit for our purposes in the current context.

\begin{df}Let $s:[I]^{<\infty}\to\mathbb{R},\ J\subset I$. We call $s$ \textbf{l-continuous} on the pair $(J,I)$ if $\forall\epsilon>0$ $\forall K\in[I]^{<\infty}$ there exists $L\in[J]^{<\infty}$ such that $d_H(L,J)\leq 2d_H(K,J)$ and $|s(K)-s(L)|<\epsilon$.
\end{df}

\begin{prp}Let $J$ be a dense subspace of $I$. If $\forall\epsilon>0\ \exists\delta>0$ such that for all $L\in[J]^{<\infty}$ and $K\in[I]^{<\infty}$,  $d_H(K,L)<\delta$ implies that $|s(K)-s(L)|<\epsilon$ (i.e. $s$ is uniformly continuous on the points of $[J]^{<\infty}$ regarding the Hausdorff pseudo-metric on $[I]^{<\infty}$) then $s$ is l-continuous on the pair $(J,I)$.
\end{prp}
\begin{proof}Choose $\delta$ for $\frac{\epsilon}{2}$ according to the definition of the uniform continuity of $s$. Now if $K\in[I]^{<\infty}$ and $\delta_1=\min\{d_H(K,J),\delta\}$ then for each $k\in K$ there exists $p_k\in J$ such that $d(k,p_k)<\delta_1$. Let $L=\{p_k:k\in K\}\in[J]^{<\infty}$. Clearly $d_H(K,L)<\delta_1$ which gives that  $d_H(L,J)\leq d_H(L,K)+d_H(K,J)<2d_H(K,J)$ and we also get that $|s(K)-s(L)|<\epsilon$.
\end{proof}

\begin{thm}\label{tlcont2ext}Let $J$ be a dense subspace of $I$. Let $s:[I]^{<\infty}\to\mathbb{R}$ and $s$ be l-continuous on the pair $(J,I)$. If ext$(s,J)$ exists then so does ext$(s,J,I)$ and they are equal.
\end{thm}
\begin{proof}Let ext$(s,J)=A$. Choose $\delta_1$ for $\frac{\epsilon}{2}$ according to the definition of ext$(s,J)$ and choose $\delta_2$ for $\frac{\epsilon}{2}$ according to the definition of the l-continuity of $s$. Let $\delta=\frac{1}{2}\min\{\delta_1,\delta_2\}$. Now if $K\in[I]^{<\infty}$ and $d_H(K,J)<\delta$ then there exists $L\in[J]^{<\infty}$ such that $d_H(L,J)\leq 2d_H(K,J)=2\delta<\delta_1$ and $|s(K)-s(L)|<\frac{\epsilon}{2}$. Which gives that $|s(K)-A|<|s(K)-s(L)|+|s(L)-A|<\frac{\epsilon}{2}+\frac{\epsilon}{2}=\epsilon$.
\end{proof}

\begin{prp}Let $J$ dense subspace of $I$ and $s$ be d-continuous. Then $s$ is l-continuous on the pair $(J,I)$ as well.
\end{prp}
\begin{proof}Let $\epsilon>0,\ K\in[I]^{n},\ d_H(K,J)=\delta_0$. For $\epsilon$ choose $\delta_1$ according to the definition of d-continuity.
Let 
\[\delta<\frac{1}{2}\min\{d(x,y):x,y\in K\}\text{ and }\delta<\min\left\{\delta_0,\delta_1\right\}.\] 
Then for each $k\in K$ there exists $p_k\in J$ such that $d(k,p_k)<\delta$. Let $L=\{p_k:k\in K\}\in[J]^{<\infty}$. Clearly $|L|=|K|$ and $d_H(L,J)\leq 2\delta_0$. Then $d_H(K,L)<\delta_1$ gives that $|s(K)-s(L)|<\epsilon$.
\end{proof}

\begin{cor}Let $J$ dense subspace of $I$. Then finite summation and finite average both l-continuous on the pair $(J,I)$.\qed
\end{cor}

%closed

Now we formulate a statement for closer of the base set.

\begin{thm}Let $J\subset I$. Let $s:[cl(J)]^{<\infty}\to\mathbb{R}$ be l-continuous  on the pair $(J,cl(J))$. If ext$(s,J)$ exists then so does and ext$(s,cl(J))$ and they are equal.
\end{thm}
\begin{proof}\ref{tlcont2ext} and \ref{pdensesp2exteq}.
\end{proof}

%--------------------------------
\subsection{Some various topics}

We are going to generalize \ref{tspacesextf}.

\begin{df}Let $J\subset I$. $J$ is called a null-subspace of $I$ regarding $s$ if $K\in[I]^{<\infty}$ then $s(K)=s(K-J)$.
\end{df}

\begin{prp}\label{pnullsspace}If $J\subset I$ is a null-subspace of $I$ regarding $s$ and $\mathrm{ext}(s,I)$ exists then so does $\mathrm{ext}(s,I-J)$ and they are equal.
\end{prp}
\begin{proof}For $\epsilon$ choose $\delta$ according to the definition of $\mathrm{ext}(s,I)$. If $K\in[I-J]^{<\infty},\ d_H(K,I-J)<\delta$ then find $L\in[J]^{<\infty}$ such that $d_H(K\cup L,I)<\delta$. Then $s(K)=s(K\cup L)$ proves the claim.
\end{proof}

\begin{df}Let $s:[I]^{<\infty}\to\mathbb{R}$ and $(I_n)$ be an ascending sequence of subsets of $I$ such that $\bigcup\limits_{n=1}^{\infty}I_n=I$. Then $(I_n)$ is called a b-converging sequence if $\lim\limits_{n\to\infty}\mathrm{ext}(s||I_n,I_n)=\mathrm{ext}(s,I)$ and $\forall\epsilon>0\ \exists\delta>0\ \exists N$ such that $n>N, K\in[I]^{<\infty}, K_n\in[I_n]^{<\infty}, d_H(K_n,I_n)<\delta, I_n\subset S_I(K,\delta)$ implies that $|s(K)-s(K_n)|<\epsilon$. 
\end{df}

\begin{prp}Let $s:[I]^{<\infty}\to\mathbb{R}$ and $(I_n)$ be an ascending sequence of subsets of $I$ such that $\bigcup\limits_{n=1}^{\infty}I_n=I$. Let $\lim\limits_{n\to\infty}\mathrm{ext}(s||I_n,I_n)=\mathrm{ext}(s,I)$. Then $(I_n)$ is a b-converging sequence iff  $\forall\epsilon>0\ \exists\delta>0\ \exists N$ such that $n>N, K\in[I]^{<\infty}, I_n\subset S_I(K,\delta)$ implies that $|s(K)-\mathrm{ext}(s,I_n)|<\epsilon$.\qed 
\end{prp}

\begin{thm}Let $I_0,I_1$ be disjoint subspaces of $I,\ I=I_0\cup I_1$ and let $s_0=s||_{I_0},\ s_1=s||_{I_1}$. Let $\tilde{I}_i=I-cl(I_{1-i})\ (i=0,1)$. Let $\hat{I}_i=cl(I_{1-i})\cap I_i\ (i=0,1)$. Let $J_{i,n}=I-S(I_{1-i},\frac{1}{n})\ (i=0,1;\ n\in\mathbb{N})$. 
Let $(J_{i,n})$ be a b-converging sequence in $\tilde{I}_i\ (i=0,1)$. Let $\hat{I}_i$ be a null-subspace of $I_i$ regarding $s_i\ (i=0,1)$.
Suppose there is a continuous function $F:\mathbb{R}\times\mathbb{R}\to\mathbb{R}$ such that $s(K)=F(s_0(K\cap I_0),s_1(K\cap I_1))$ ($K\in [I]^{<\infty}$). If ext$(s_0,I_0)$ and ext$(s_1,I_1)$ exist then so does ext$(s,I)$ and ext$(s,I)=F($ext$(s,I_0),$ext$(s,I_1))$.  
\end{thm}
\begin{proof}Let $\tilde{s}_i=s_i||_{\tilde{I}_i}\ (i=0,1)$ and $\tilde{s}_{i,n}=s_i||_{J_{i,n}}\ (i=0,1;\ n\in\mathbb{N}))$. Observe that $I_i=\tilde{I}_i\cup^*\hat{I}_i\ (i=0,1)$

Let ext$(s_i,I_i)=A_i\ (i=0,1)$. Let $\epsilon>0$ and choose $\epsilon_1$ according to the continuity of $F$ at $(A_0,A_1)$: if $x_i\in S(A_i,\epsilon_1)\ (i=0,1)$ then $|F(x_0,x_1)-F(A_0,A_1)|<\epsilon$. 

Note that $(J_{i,n})$ is ascending and $\bigcup\limits_{n=1}^{\infty}J_{i,n}=\tilde{I}_i\ (i=0,1)$. Hence there is $N$ such that $n\geq N$ implies that $\big|\mathrm{ext}(\tilde{s}_{i,n},J_{i,n})-\mathrm{ext}(\tilde{s}_i,\tilde{I}_i)\big|<\frac{\epsilon_1}{2}\ (i=0,1)$. Choose $\delta_i<\frac{1}{N}$ such that $K_i\in[\tilde{I}_i]^{<\infty},\ J_{i,N}\subset S_{\tilde{I}_i}(K_i,\delta_i)$ implies that $\big|\tilde{s}_i(K_i)-\mathrm{ext}(\tilde{s}_{i,N},J_{i,N})\big|<\frac{\epsilon_1}{2}\ (i=0,1)$.

Let $\delta=\min\{\delta_0,\delta_1\}$. If $K\in[I]^{<\infty}$ such that $d_H(K,I)<\delta$ then clearly $J_{i,N}\subset S_{\tilde{I}_i}(K\cap \tilde{I}_i,\delta)\ (i=0,1)$. Then
\[|s_i(K\cap I_i)-A_i|\leq|s_i(K\cap I_i)-s_i(K\cap \tilde{I}_i)|
+|s_i(K\cap \tilde{I}_i)-\mathrm{ext}(\tilde{s}_{i,N},J_{i,N})|+
\]
\[
+|\mathrm{ext}(\tilde{s}_{i,N},J_{i,N})-\mathrm{ext}(\tilde{s}_i,\tilde{I}_i)|
+|\mathrm{ext}(\tilde{s}_i,\tilde{I}_i)-A_i|\leq 0+\frac{\epsilon_1}{2}+\frac{\epsilon_1}{2}+0=\epsilon_1\]
because $\hat{I}_i$ be a null-subspace of $I_i$ and $\mathrm{ext}(\tilde{s}_i,\tilde{I}_i)=A_i\ (i=0,1)$ by \ref{pnullsspace}. Therefore $|s(K)-F(A_0,A_1)|<\epsilon$.
\end{proof}

%\begin{rem}Actually we used less than the b-continuity of $I$: We just used that $\mathrm{ext}(\tilde{s}_i,\tilde{I}_i)$ is an accumulation point of $(\mathrm{ext}(\tilde{s}_i,J_{i,n},\tilde{I}_i))$.
%\end{rem}

We now investigate the relation of extensions and continuity.

\begin{lem}\label{luc1}Let $f:J\to I$ be a uniformly continuous surjective function and $\epsilon>0$. Let $\delta$ be chosen according to the uniform continuity of $f$. If $K\in[J]^{<\infty}$ and $\delta$-sdense then $f(K)$ is $\epsilon$-sdense.
\end{lem}
\begin{proof}Assume the contrary. Then there is $x\in I$ such that $x\notin S(f(K),\epsilon)$. Choose $y\in f^{-1}(\{x\})$ and then there is $z\in K$ such that $d(z,y)<\delta$. But we get that $d(f(z),f(y))<\epsilon$ which is a contradiction.
\end{proof}

\begin{thm}Let $s_i:[I]^{<\infty}\to\mathbb{R}$, $f:J\to I$ be a uniformly continuous surjective function and let $s_j=s_i\circ f$. Then if ext$(s_i,I)$ exists then so does ext$(s_j,J)$ and they are equal.
\end{thm}
\begin{proof}Let ext$(s_i,I)=A$ and $\epsilon>0$. Choose $\delta$ according to the definition of ext$(s_i,I)$. Then for $\delta$ choose $\xi$ according to the definition of uniform continuity. If $K\in[J]^{<\infty}$ and $\xi$-sdense then $f(K)$ is $\delta$-sdense by \ref{luc1}. Clearly $s_j(K)=s_i(f(K))$ hence $|s_j(K)-A|=|s_i(f(K))-A|<\epsilon$.
\end{proof}

Now we investigate some connection between extension and products.

\begin{thm}Let $\langle I_1,d_1\rangle,\ \langle I_2,d_2\rangle$ be totally bounded pseudo-metric spaces, $I=I_1\times I_2,\ d=d_1+d_2$. Let $F:\mathbb{R}\times\mathbb{R}\to\mathbb{R}$ be continuous. Let $s_i:[I_i]^{<\infty}\to\mathbb{R}\ (i=1,2),\ s:[I]^{<\infty}\to\mathbb{R}$ such that if $K_i\in[I_i]^{<\infty}\ (i=1,2)$ then let $s(K_1\times K_2)=F(s_1(K_1),s_2(K_2))$. Now if ext$(s_i,I_i)\ (i=1,2)$ and ext$(s,I)$ all exist then ext$(s,I)=F\left(\mathrm{ext}(s_1,I_1),\mathrm{ext}(s_2,I_2)\right)$.
\end{thm}
\begin{proof}Let ext$(s_i,I_i)=A_i$, ext$(s,I)=A$ and $\epsilon>0$. Find $\epsilon_1$ such that if $|A_i-x_i|<\epsilon_1\ (i=1,2)$ then $|F(A_1,A_2)-F(x_1,x_2)|<\frac{\epsilon}{2}$. Choose then $\delta$ such that $K_i\in[I_i]^{<\infty},\ d_H(K_i,I_i)<\delta$ implies that $|s_i(K_i)-A_i|<\epsilon_1\ (i=1,2)$ and $K\in[I]^{<\infty},\ d_H(K,I)<\delta$ implies that $|s(K)-A|<\frac{\epsilon}{2}$. Now if $K_i\in[I_i]^{<\infty},\ d_H(K_i,I_i)<\frac{\delta}{2}\ (i=1,2)$ then $d_H(K_1\times K_2,I_1\times I_2)<\delta$. This gives that $|s(K_1\times K_2)-F(A_1,A_2)|=|F(s_1(K_1),s_2(K_2))-F(A_1,A_2)|<\frac{\epsilon}{2}$ and $|s(K_1\times K_2)-A|<\frac{\epsilon}{2}$ hence $|F(A_1,A_2)-A|<\epsilon$.
\end{proof}

%-------------------------------------------------------------------------------------------------------------------------------------------Some generalizations----------------------------------
\section{Some generalizations}

In \ref{pintinmeassp} we had the restriction that $\mu(X)<\infty$. We are now going to get rid of this.

First let us allow that the function $s$ can also take the value $+\infty$ in the definitions of ext \ref{dext} and \ref{dextinfty} i.e. $s:[I]^{<\infty}\to\mathbb{R}\cup\{+\infty\}$. Then modify slightly the definition of a function being d-increasing .

\begin{df}\label{ddinc}A function $s:[I]^{<\infty}\to\mathbb{R}\cup\{+\infty\}$ is called d-increasing if $\forall\epsilon>0$ and for all $K\in[I]^{<\infty}$ there is $\delta>0$ such that $L\in[I]^{<\infty}$ being $\delta$-sdense implies that $s(K)-\epsilon\leq s(L)$.
\end{df}

\begin{rem}It can be readily seen that \ref{pdincsup} remains valid for these extended versions of base function $s$ and notion d-increasing.
\end{rem}

\begin{prp}\label{pintinmeassp3}Let $\langle X,\mu\rangle$ be a measure space. Let $f:X\to\mathbb{R}$ be a measurable function such that $0\leq f<M\in\mathbb{R}$. Let $I=(0,M)$ and $s:[I]^{<\infty}\to\mathbb{R}\cup\{+\infty\}$ defined as in \ref{pintinmeassp}. Then $s$ is d-increasing.
\end{prp}
\begin{proof}We have to manage when $\mu(X)=+\infty$.
Let $H\in[I]^{<\infty},\ \epsilon>0$. Then there are two cases.
\begin{enumerate}
\item Let $s(H)<\infty$. Then clearly $\mu\Big(f^{-1}\big([\min H,M)\big)\Big)<\infty$. Let 
\[\delta<\frac{\epsilon}{\mu\Big(f^{-1}\big([\min H,M)\big)\Big)}.\]
Let $K\in[I]^{<\infty},\ d_H(K,I)<\delta$. If $s(K)=+\infty$ then we are done. 
If $s(K)<\infty$ then replace $X$ with $f^{-1}\big([\min H,M)\big)$ in the proof of \ref{pintinmeassp} and use the trivial fact that $c_i\in H$ implies that $c_i\geq \min H$.

\item Let $s(H)=+\infty$. Let $\delta=\frac{\min H}{2}$ and $K\in[I]^{<\infty},\ d_H(K,I)<\delta$. Then $\mu\Big(f^{-1}\big([\min H,M)\big)\Big)=\infty$ and $[\min H,M)\subset[\min K,M)$ gives that $\mu\Big(f^{-1}\big([\min K,M)\big)\Big)=\infty$ therefore $s(K)=+\infty$ as well.\qedhere
\end{enumerate}
\end{proof}

\begin{thm}Using the notations of \ref{pintinmeassp3} ext$(s,I)$ exists and equals to $\int\limits_X f d\mu$.\qed
\end{thm}

In order to manage the integral of unbounded functions on measure spaces we need to generalize the notion of $\delta$-sdense sets. It is needed because $(0,+\infty)$ is not totally bounded.

\begin{df}Let $0<\delta<1$. A subset $K$ of $(0,+\infty)$ is called $\delta$-dense in $(0,+\infty)$ if $K\cap\left(0,\frac{1}{\delta}\right)$ is $\delta$-dense in $\left(0,\frac{1}{\delta}\right)$. $K$ is $\delta$-sdense in $(0,+\infty)$ if there is $\delta'<\delta$ such that $K$ is $\delta'$-dense in $(0,+\infty)$. 
\end{df}

\begin{prp}If $K\subset (0,+\infty)$ is $\delta$-dense in $(0,+\infty)$ then both sets $K\cap(0,1)$ and $\tilde K=\left\{\frac{1}{k}:k\in K\cap[1,+\infty)\right\}$ are $2\delta$-dense in $(0,1)$.
\end{prp}
\begin{proof}Let $a=\max K\cap(0,1),\ b=\min K\cap[1,+\infty),\ c=\max K\cap[1,+\infty)$.

Regarding $K\cap(0,1)$ $S(b,\delta)$ can cover some part of $(0,1)$. But $(0,1)-S(b,\delta)\subset S(K\cap(0,1),\delta)$ which gives that $(0,1)\cap S(b,\delta)\subset S(a,2\delta)$. Hence $(0,1)\subset S(K\cap(0,1),2\delta)$.

For $\tilde K$ note first that if $x,y\geq 1,\ |x-y|<\delta$ then $\left|\frac{1}{x}-\frac{1}{y}\right|<\delta$ too. Then $\frac{1}{\delta}\in S(c,\delta)$ which implies that $0\in S(\frac{1}{c},2\delta)$. Moreover similar reasoning as in the first part of the proof gives that $S(a,\delta)\cap [1,+\infty)\subset S(b,2\delta)$ hence $1\in S(\frac{1}{b},2\delta)$. Therefore $(0,1)\subset S(\tilde K,2\delta)$.
\end{proof}

\begin{prp}Let $K\subset (0,+\infty)$ and $0<\delta<1$. If both $K\cap(0,1)$ and $\tilde K=\left\{\frac{1}{k}:k\in K\cap[1,+\infty)\right\}$ are $\delta^3$-dense in $(0,1)$ then $K$ is $\delta$-dense is $(0,+\infty)$.
\end{prp}
\begin{proof}Obviously $K\cap(0,1)$ is $\delta$-dense since $\delta^3<\delta$.

Regarding $K\cap[1,+\infty)$ we have to prove that $K\cap[1,+\infty)$ is $\delta$-dense in $\left[1,\frac{1}{\delta}\right)$. Therefore it is enough to show the following: 

If $x,y\in\tilde K,\ x<y,\ \frac{1}{x}<\frac{1}{\delta},\ \frac{1}{y}<\frac{1}{\delta},\ y-x<2\delta^3$ then $\frac{1}{x}-\frac{1}{y}<2\delta$. 

But $\frac{1}{x}<\frac{1}{\delta},\frac{1}{y}<\frac{1}{\delta}$ implies that $\frac{1}{xy}<\frac{1}{\delta^2}$ hence we get that
\[\frac{1}{x}-\frac{1}{y}=\frac{y-x}{xy}<\frac{2\delta^3}{\delta^2}<2\delta.\qedhere\]
\end{proof}

\begin{df}The definition of ext$(s,I)$ (see \ref{dext}) can be slightly extended by allowing $I=(0,+\infty)$ and replacing $d_H(K,I)<\delta$ with $K$ being $\delta$-sdense in $(0,+\infty)$.
\end{df}

\begin{df}The definition of a function being d-increasing (d-decreasing) (see \ref{ddinc}) can be slightly extended by allowing $I=(0,+\infty)$.
\end{df}

\begin{prp}\label{pintinmeassp2}Let $\langle X,\mu\rangle$ be a measure space such that $\mu(X)<\infty$. Let $f:X\to\mathbb{R}$ be a measurable function such that $0\leq f$. Let $I=(0,+\infty)$ and $s:[I]^{<\infty}\to\mathbb{R}$ be defined as follows. If $H\in[I]^{<\infty},\ H=\{a_1,\dots,a_n\},\ a_1<\dots<a_n$ then let $a_0=0,a_{n+1}=+\infty$ and set
\[s(H)=\sum\limits_{i=0}^{n}a_i\mu\Big(f^{-1}\big([a_i,a_{i+1})\big)\Big).\]
Then $s$ is d-increasing.
\end{prp}
\begin{proof}The proof of proposition \ref{pintinmeassp} can be copied with the following two modifications:

1. $\delta$ has to be chosen as $\delta<\min\{\frac{1}{\max H},\frac{\epsilon}{\mu(X)}\}$

2. $d_H(K,I)<\delta$ has to be replaced with $K$ being $\delta$-sdense on $(0,+\infty)$.
\end{proof}

\begin{thm}Let $\langle X,\mu\rangle$ be a measure space such that $\mu(X)<\infty$. Let $f:X\to\mathbb{R}$ be a measurable function such that $0\leq f$. Let $I=(0,+\infty)$ and $s:[I]^{<\infty}\to\mathbb{R}$ be defined as in \ref{pintinmeassp2}. Then ext$(s,I)$ exists and equals to $\int\limits_X f d\mu$.
\end{thm}
\begin{proof}
By \ref{pdincsup} ext$(s,I)$ exists and equals to $\sup\{s(K):K\in [I]^{<\infty}\}$ which is known to equal to $\int\limits_X f d\mu$.
\end{proof}
%---------------------------------------------------------------------------------------------------------------------------------------------------------the bibliography -----------------

{\footnotesize

%------------------------------------------------------------------------------------------------------------------------------------------------------------address---------------
\noindent
Dennis G\'abor College, Hungary 1119 Budapest Fej\'er Lip\'ot u. 70.

\noindent E-mail: losonczi@gdf.hu, alosonczi1@gmail.com\\
}
\end{document}